\documentclass[11pt,a4paper]{article}
\usepackage[T1]{fontenc}
\usepackage{times}
\usepackage{amsmath,amssymb,amsthm}
\usepackage{txfonts}
\usepackage[margin=28mm]{geometry}
\usepackage{microtype}
\usepackage{booktabs,array,longtable}
\usepackage{enumitem}
\usepackage{needspace}
\usepackage{indentfirst}

\usepackage[hidelinks]{hyperref}
\usepackage{xcolor}
\hypersetup{pdftitle={k-trees with second largest eigenvalue at most one},
 pdfauthor={Zihao Geng, Jongyook Park, Ji-Ming Guo},
 pdfsubject={Structural classification and the second largest eigenvalues nearest to one},
 pdfkeywords={k-tree, second largest eigenvalue, forbidden induced subgraph, multiplicity}}
\setlist[enumerate]{label=(\roman*),leftmargin=*,itemsep=2pt,topsep=4pt}
\newtheorem{theorem}{Theorem}
\newtheorem{lemma}[theorem]{Lemma}
\newtheorem{proposition}[theorem]{Proposition}
\newtheorem{corollary}[theorem]{Corollary}
\theoremstyle{remark}
\newtheorem{remark}[theorem]{Remark}
\newcommand{\ind}{\leq_{\mathrm{ind}}}
\newcommand{\aone}{\boldsymbol{a}}
\newcommand{\Fcal}{\mathcal{F}}
\newcommand{\R}{\mathbb{R}}
\allowdisplaybreaks[1]
\renewenvironment{abstract}{%
 \par\addvspace{\baselineskip}\noindent
 {\large\bfseries\abstractname\par}%
 \nobreak\noindent\ignorespaces
}{\par\vspace{\baselineskip}\vspace{0.1cm}}
\title{$k$-trees with second largest eigenvalue at most one}
\author{%
 {\Large Zihao Geng$^{a,b}$,\,\,
 Ji-Ming Guo$^{a}$,\,\,
 Jongyook Park$^{b}$%
 \thanks{Corresponding author.}}\\[2mm]
 \parbox[t]{0.96\textwidth}{%
 \centering\fontsize{10}{12}\selectfont\itshape
 \setlength{\parindent}{0pt}\setlength{\parskip}{4pt}
 \textsuperscript{a}\,Department of Mathematics,
 East China University of Science and Technology, Shanghai 200237, China\endgraf
 \textsuperscript{b}\,Department of Mathematics,
 Kyungpook National University, Daegu 41566, Republic of Korea}}
\date{}
\begin{document}
\maketitle
\begingroup
\makeatletter
\let\@footnotetext\H@@footnotetext
\makeatother
\renewcommand{\thefootnote}{}
\urlstyle{same}
\footnotetext{\textit{E-mail addresses:}\quad
 \href{mailto:gzh635516986@163.com}{\nolinkurl{gzh635516986@163.com}} (Z.H. Geng),
 \href{mailto:jimingguo@hotmail.com}{\nolinkurl{jimingguo@hotmail.com}} (J.-M. Guo),
 \href{mailto:jongyook@knu.ac.kr}{\nolinkurl{jongyook@knu.ac.kr}} (J. Park).}
\endgroup
\begin{abstract}
The problem of characterizing  graphs whose second largest adjacency eigenvalue is at most $1$ was posed by Hoffman and first studied by Cvetkovi\'c in 1982.
In this paper,  all $k$-trees with second largest eigenvalue at most $1$ are characterized.
We determine all equality cases and the multiplicity of the eigenvalue $1$.
For each fixed $k\geq2$, we also determine the $k$-trees attaining the largest  second largest eigenvalue below $1$ and the smallest above $1$.
Finally, we determine all minimal forbidden induced subgraphs among $k$-trees with $k\geq 1$, each of which has second largest eigenvalue greater than $1$, but deleting any vertex gives a graph whose second largest eigenvalue is at most $1$.
\end{abstract}
\par
\vspace{1.5\baselineskip}


\section{Introduction}\label{sec:introduction}
In this paper, we study $k$-trees whose second largest adjacency eigenvalue is at most $1$. All graphs considered here are finite, undirected and simple. For a graph $G$ of order $n\geq2$, let $V(G)$ and $E(G)$ denote its vertex set and edge set, respectively, and let $A(G)$ be its \emph{adjacency matrix}. We call the eigenvalues of $A(G)$ the \emph{(adjacency) eigenvalues} of $G$ and write them as $\lambda_1(G)\geq\lambda_2(G)\geq\cdots\geq\lambda_n(G)$. In particular, $\lambda_2:=\lambda_2(G)$ is the \emph{second largest eigenvalue} of $G$.

We use $K_n$ and $P_n$ for the complete graph and the path on $n$ vertices, respectively, and $\overline{K_n}$ for the graph with $n$ vertices and no edges, where $K_0$ denotes the graph with no vertices. A set of $n$ pairwise adjacent vertices is called an \emph{$n$-clique}. The \emph{join} $G\vee H$ of two vertex-disjoint graphs is obtained from their disjoint union by adding every edge between $V(G)$ and $V(H)$. We denote their disjoint union by $G\cup H$. We write $H\ind G$ if $H$ is isomorphic to an \emph{induced subgraph} of $G$. For graphs $F_1,\ldots,F_s$, a graph is $(F_1,\ldots,F_s)$-free if it contains no induced subgraph isomorphic to any $F_i$. Let $K_t-e$ denote the graph obtained by deleting one edge from $K_t$. Let $T_{1,2}$ denote the tree of diameter $3$ with exactly three pendant vertices.

The problem of characterizing connected graphs with $\lambda_2\leq1$ was posed by Hoffman and first studied by Cvetkovi\'c~\cite{Cvetkovic}.
This problem has been studied for several classes of graphs.
Trees~\cite{Neumaier}, unicyclic graphs~\cite{Xu2004}, bicyclic graphs~\cite{Guo2005}, tricyclic graphs~\cite{LiYang2011}, bipartite graphs~\cite{Petrovic1991}, line graphs~\cite{PetrovicMilekic1998}, generalized line graphs~\cite{PetrovicMilekic2000}, generalized $\theta$-graphs~\cite{GaoHuang2008}, connected graphs with three distinct eigenvalues~\cite{ChengGavrilyukGreavesKoolen2016,ChengGreavesKoolen2018}, $r$-regular graphs of order $n$ with $r\leq10$ or $n\geq2r$~\cite{KoledinStanic2013,Stanic2010,WangZhu2026}, nested split graphs~\cite{MilatovicStanic2012,Stanic2009}, connected $(K_{1,3},K_5-e)$-free graphs~\cite{LiuChenStanic2021}, connected $(K_4-e)$-free graphs~\cite{LiuChenStanic2024}, connected outerplanar graphs~\cite{LiSun}, connected $T_{1,2}$-free planar graphs~\cite{GengGuoWang2024}, connected triangle-free planar graphs~\cite{ChengLi2025}, connected $K_4$-minor-free graphs~\cite{Cheng} and connected planar graphs~\cite{GengWangLiLiChenGuo2026} with $\lambda_2\leq1$ have been characterized.
Non-bipartite graphs of girth at least $4$ have also been studied with $\lambda_2<1$~\cite{XuShao2006}.
In addition, the $2$-partially distance-regular graphs whose local graphs all have second largest eigenvalue at most $1$ have been characterized~\cite{ZhangLiangKoolen2022}.

For a positive integer $k$, a \emph{$k$-tree} is obtained from $K_k$ by repeatedly adding a vertex adjacent precisely to the vertices of an existing $k$-clique. 
We include the initial graph $K_k$ in this definition. The complete split graph $K_k\vee\overline{K_t}$ is a $k$-tree for every integer $t\geq0$, and every maximal outerplanar graph of order at least three is a $2$-tree.
The structure of $k$-trees was studied by Beineke and Pippert~\cite{BeinekePippert} and Rose~\cite{Rose}.
Since a $1$-tree is a tree and a $2$-tree is $K_4$-minor-free,  the case $k\leq2$ of our classification is contained in ~\cite{Cheng,Neumaier}.
The main new range in our classification is therefore $k\geq3$, while the argument also covers $k=2$. A $k$-tree of order greater than $k$ has clique number $k+1$, so the result applies to graphs of arbitrarily large clique number. 
Since complete graphs of order at least two have second largest eigenvalue $-1$, the presence of large cliques does not itself prevent the bound $\lambda_2\leq1$. 
Our classification describes which attachments to these cliques preserve the bound. 
We also determine the possible values of $\lambda_2$ nearest to $1$, rather than only the graphs satisfying the single bound $\lambda_2\leq1$.

Let $Q=\{q_1,\ldots,q_{k+1}\}$ induce $K_{k+1}$. For nonnegative integers $a_1,\ldots,a_{k+1}$, take pairwise disjoint sets $X_1,\ldots,X_{k+1}$, disjoint from $Q$, with $|X_i|=a_i$. Require $X_1\cup\cdots\cup X_{k+1}$ to be an \emph{independent set}. Join every vertex of $X_i$ to every vertex of $Q\setminus\{q_i\}$, and add no other edges. Denote the resulting $k$-tree by $B_k(a_1,\ldots,a_{k+1})$. We call $Q$ its \emph{core} and the sets $X_i$ its \emph{attachment classes}. For integers $0\leq s\leq r$, let $S_{r,s}$ be the tree obtained from the star $K_{1,r}$ by subdividing exactly $s$ pendant edges once. Thus $S_{r,0}=K_{1,r}$. In particular, $S_{0,0}=K_1$. For $k\geq2$, put $r_k=\left\lfloor7+\frac{2}{k-1}\right\rfloor$, and define
\begin{equation}\label{eq:hosts}
 H_{k,1}=B_k(2,\ldots,2),  \qquad H_{k,2}=B_k(3,1,1,0,\ldots,0),  \qquad H_{k,3}=K_{k-1}\vee S_{r_k,r_k}.
\end{equation}
In particular, $r_2=9$, $r_3=8$, and $r_k=7$ for $k\geq4$. {The following theorem gives the classification in terms of these three graphs and the graphs $K_k\vee\overline{K_t}$ with $t\geq0$.}

\begin{theorem}\label{thm:classification}
Let $k\geq2$, and let $G$ be a $k$-tree. Then $\lambda_2(G)\leq1$ if and only if one of the following holds:
\begin{enumerate}
\item $G\cong K_k\vee\overline{K_t}$ for some integer $t\geq0$;
\item $G\ind H_{k,j}$ for some $j\in\{1,2,3\}$.
\end{enumerate}
\end{theorem}

The proof uses interlacing to exclude a finite collection of induced subgraphs and then reduces $G$ to one of the families $B_k(a_1,\ldots,a_{k+1})$ and $K_{k-1}\vee S_{r,s}$. Characteristic polynomials verify the required eigenvalue bounds, and the remaining parameter restrictions follow from the graph structure. The forbidden graphs are themselves $k$-trees, and we show that their list is complete and minimal. Then eigenvalue equations determine all equality cases and the exact multiplicity of $1$.

Our second main result determines the possible values of $\lambda_2$ immediately below and above $1$. Put $\sigma_2=7$ and $\sigma_k=6$ for $k\geq3$, and define $L_k=B_k(\sigma_k,1,0,\ldots,0)$ and $E_k=B_k(r_k,1,0,\ldots,0)$.

\begin{theorem}\label{thm:gap}
For each $k\geq2$, the quantities
\[
 \begin{aligned}
 \beta_k&=\max\{\lambda_2(G):G\text{ is a $k$-tree and }\lambda_2(G)<1\},\\
 \gamma_k&=\min\{\lambda_2(G):G\text{ is a $k$-tree and }\lambda_2(G)>1\}
 \end{aligned}
\]
exist and satisfy $\beta_k=\lambda_2(L_k)$ and $\gamma_k=\lambda_2(E_k)$. Up to isomorphism, $L_k$ is the unique $k$-tree whose second largest eigenvalue is $\beta_k$, and $E_k$ is the unique $k$-tree whose second largest eigenvalue is $\gamma_k$.
\end{theorem}

Thus $\lambda_2(G)\notin(\beta_k,1)\cup(1,\gamma_k)$ for every $k$-tree $G$. This determines the nearest possible values on both sides of $1$, together with the unique graphs attaining them. We will also show that $\beta_k$ tends to $1$ and that the product of $k$ and the difference between $1$ and $\beta_k$ tends to $2/35$. In contrast, $\gamma_k$ tends to a number strictly greater than $1$. In fact, $\gamma_4$ is the minimum of $\lambda_2(G)$ over all $k$-trees $G$ with $k\geq2$ and $\lambda_2(G)>1$.

{For each fixed $k$, Theorem~\ref{thm:minimal} gives a list of forbidden induced subgraphs that are themselves $k$-trees. Their orders depend on $k$. In Section~\ref{sec:uniform}, we obtain one irredundant list of nine graphs that applies to all $k\geq1$. These are also precisely the $k$-trees with $\lambda_2>1$ for which every proper induced subgraph of order at least two has second largest eigenvalue at most $1$, without requiring the subgraph to be a $k$-tree. Consequently, every $k$-tree with $\lambda_2>1$ contains a connected induced subgraph of order at most $13$ satisfying the same inequality, and $13$ is best possible. We use Neumaier's tree classification for the case $k=1$. Its elementary spectral consequences also explain why the largest value below $1$ in Theorem~\ref{thm:gap} requires $k\geq2$.}

The paper is organized as follows. In Section~\ref{sec:spectral}, we give the characteristic polynomials and interlacing results needed in the proofs. Section~\ref{sec:structure} proves Theorem~\ref{thm:classification} and determines the minimal forbidden $k$-trees. In Section~\ref{sec:boundary}, we determine the equality cases, the multiplicity of $1$, {and the sharp order bound for graphs not isomorphic to $K_k\vee\overline{K_t}$ for any $t\geq0$.} Section~\ref{sec:gap} proves Theorem~\ref{thm:gap} and describes the limiting behavior of $\beta_k$ and $\gamma_k$. {Finally, Section~\ref{sec:uniform} recalls the necessary tree results and gives the common forbidden induced subgraphs for all $k\geq1$, their minimality under arbitrary vertex deletion, and the sharp bound of $13$ on their orders.}

\section{Preliminaries}\label{sec:spectral}
In this section, we give the polynomial factorizations used to verify the classification and to compare the values of $\lambda_2$ nearest to $1$.
The \emph{neighborhood} of a vertex $v$ in a graph $G$ is denoted by $N_G(v)$, and $G[U]$ is the subgraph induced by $U\subseteq V(G)$. We write $\phi(G;x)=\det(xI-A(G))$ for the \emph{characteristic polynomial} of $G$. For $\theta\in\R$, let $m_G(\theta)=\dim\ker(A(G)-\theta I)$ denote its \emph{multiplicity} as an adjacency eigenvalue, with $m_G(\theta)=0$ when $\theta$ is not an eigenvalue.

We use interlacing in the following form: if $H\ind G$ and $|V(H)|\geq2$, then $\lambda_2(H)\leq\lambda_2(G)$. We also use $\lambda_1(H)\leq\lambda_1(G)$. Interlacing further gives $\lambda_2(G)\leq\lambda_1(G-v)$ for every vertex $v$ of a graph $G$ of order at least two. We use the equality condition as well. If $\lambda_i(H)=\lambda_i(G)$, then there is a nonzero eigenvector of $H$ for this eigenvalue whose extension by zero on $V(G)\setminus V(H)$ is an eigenvector of $G$, see~\cite[Theorem~2.5.1]{BrouwerHaemers}.

In particular, suppose that $\theta=\lambda_2(H)$ is \emph{simple}, that is, $m_H(\theta)=1$, and let $\boldsymbol z$ be a corresponding nonzero eigenvector. Let $H'$ be obtained from $H$ by adding a vertex $v$ with neighborhood $W\subseteq V(H)$. If $\sum_{u\in W}z_u\ne0$, then $\lambda_2(H')>\theta$. To see this, assume that equality holds in the interlacing inequality $\lambda_2(H')\geq\lambda_2(H)$. Then, there exists a nonzero eigenvector $\boldsymbol y$ of $H$ corresponding to $\theta$ whose extension by zero at $v$ is an eigenvector of $H'$. Since the eigenspace of $H$ corresponding to $\theta$ has dimension one, the two nonzero eigenvectors $\boldsymbol y$ and $\boldsymbol z$ must satisfy $\boldsymbol y=c\boldsymbol z$ for some scalar $c\ne0$. The eigenvalue equation at the new vertex is therefore $0=\sum_{u\in W}y_u=c\sum_{u\in W}z_u$, a contradiction.

Recall that a partition $V(G)=V_1\cup\cdots\cup V_m$ is \emph{equitable} if every vertex of $V_i$ has the same number $q_{ij}$ of neighbors in $V_j$ $(1\leq i,j\leq m)$. The matrix $(q_{ij})$ is its \emph{quotient matrix}. Every eigenvalue of this matrix is an eigenvalue of $G$, and an eigenvector is obtained by assigning the $i$th quotient coordinate to every vertex of $V_i$. The quotient matrix accounts for $m$ eigenvalues, counted with multiplicity. The other eigenvalues have eigenvectors whose coordinates sum to zero on each part,  see~\cite[Section~2.3]{BrouwerHaemers}. When we calculate a characteristic polynomial, we also identify these remaining eigenvalues, so that the factorization describes the full spectrum rather than only the quotient spectrum.

Two elementary observations will be used repeatedly. If two nonadjacent vertices have the same neighbors, the vector with coordinates $1$ and $-1$ at these vertices and zero elsewhere is an eigenvector for $0$. If two adjacent vertices have the same neighbors outside the pair, the same vector is an eigenvector for $-1$. Every other vertex is adjacent either to both selected vertices
or to neither of them. So, its neighboring coordinates sum to $1+(-1)=0$ in the first case and to $0$ in the second, i.e., in each case the contributions at every other vertex cancel. Thus, an independent set of $s$ vertices with identical neighborhoods gives $s-1$ independent eigenvectors for $0$, and a clique of $s$ vertices with identical neighbors outside the clique gives $s-1$ independent eigenvectors for $-1$. These vectors are obtained by fixing one vertex and pairing it with each of the other $s-1$ vertices.

We bound the number of positive roots of a polynomial by considering the number of sign changes in its nonzero coefficients, read in descending powers. This is known as Descartes' rule of signs. A root of multiplicity $m$ is counted $m$ times. In particular, if the rule allows at most two positive roots and a sign calculation gives roots in two disjoint positive intervals, there is exactly one root in each interval and both are simple. To locate eigenvalues greater than $1$, we apply this argument after substituting $x=y+1$. Decimal endpoints such as $0.97$ and $1.03$ stand for the exact fractions $97/100$ and $103/100$.

For calculations involving a leaf $v$ with unique neighbor $u$, we use
\begin{equation}\label{eq:leaf-recurrence}
 \phi(G;x)=x\phi(G-v;x)-\phi(G-\{u,v\};x),
\end{equation}
where $\phi(K_0;x)=1$. Here is a simple explanation. In the determinant of $xI-A(G)$, the row indexed by $v$ has only two nonzero entries, namely $x$ in column $v$ and $-1$ in column $u$. The terms choosing the first entry sum to $x\phi(G-v;x)$. If the entry in column $u$ is chosen instead, column $v$ must use the entry $-1$ in row $u$. The transposition of $u$ and $v$ gives the minus sign, and the other entries give $\phi(G-\{u,v\};x)$. These are all possible nonzero terms. We also use that the characteristic polynomial of a disjoint union is the product of those of its components, since its adjacency matrix is block diagonal.

We first calculate the characteristic polynomials of $H_{k,1}$ and $H_{k,2}$.

\begin{lemma}\label{lem:hostB}
For $k\geq2$,
\begin{align}
 \phi(H_{k,1};x)&=x^{k+1}(x-2k)(x+k)(x-1)^k(x+2)^k,\label{eq:hostone}\\
 \phi(H_{k,2};x)&=x^2(x+1)^{k-3}(x^2+x-1)p_k(x) \qquad(k\geq3),\label{eq:hosttwo}\\
 \phi(H_{2,2};x)&=x^2(x^2+x-1)(x^4-x^3-11x^2-2x+12),\label{eq:hosttwo2}
\end{align}
where
\begin{equation}\label{eq:pk}
 p_k(x)=x^5+(2-k)x^4+(2-7k)x^3+(1-7k)x^2+(8k-6)x+6k.
\end{equation}
\end{lemma}
\begin{proof}
First consider $H_{k,1}$. Each set $X_i$ has two vertices with the same neighbors, so their coordinate difference is an eigenvector for $0$, i.e., for each $i$, let $x_i,x_i'$ be the two vertices of $X_i$, then the vector with coordinates $1$ at $x_i$, $-1$ at $x_i'$, and $0$ at every other vertex is an eigenvector for $0$, since $x_i$ and $x_i'$ have the same neighbors. This gives $k+1$ independent eigenvectors for $0$.

For each $i\in\{1,\ldots,k\}$, assign coordinates $1,-1$ to $q_i,q_{k+1}$, coordinates $-1$ to the two vertices of $X_i$, coordinates $1$ to the two vertices of $X_{k+1}$, and zero elsewhere. At $q_i$ the sum of the neighboring coordinates is $-1+2=1$, and at a vertex of $X_i$ it is $-1$. The equations at $q_{k+1}$ and $X_{k+1}$ are the corresponding negatives, and at every other vertex the contributions cancel. Thus this is an eigenvector for $1$. The $k$ choices of $i$ are linearly independent, as is seen from their coordinates at $q_1,\ldots,q_k$.

Similarly, assign coordinates $2,-2$ to $q_i,q_{k+1}$, coordinates $1$ on $X_i$, coordinates $-1$ on $X_{k+1}$, and zero elsewhere. The sums of the neighboring coordinates at $q_i$ and a vertex of $X_i$ are $-4$ and $-2$, respectively. The other equations follow in the same way, so this gives $k$ independent eigenvectors for $-2$.

The partition with two parts $Q$ and $X_1\cup\cdots\cup X_{k+1}$ is equitable, with quotient matrix $\left(\begin{smallmatrix}k&2k\\k&0\end{smallmatrix}\right)$. Its characteristic polynomial is $(x-2k)(x+k)$. The eigenvectors just constructed have coordinate sum zero on each of these two parts and are independent of the two quotient eigenvectors. Therefore we have accounted for $(k+1)+k+k+2=3k+3$ eigenvalues, with multiplicity, which is the order of $H_{k,1}$. This proves~\eqref{eq:hostone}.

Now let $k\geq3$ and consider $H_{k,2}$. The partition with parts $\{q_1\}$, $\{q_2,q_3\}$, $Q\setminus\{q_1,q_2,q_3\}$, $X_1$, and $X_2\cup X_3$ is equitable. Counting neighbors in these parts gives the quotient matrix
\[
 \begin{pmatrix}
 0&2&k-2&0&2\\
 1&1&k-2&3&1\\
 1&2&k-3&3&2\\
 0&2&k-2&0&0\\
 1&1&k-2&0&0
 \end{pmatrix},
\]
whose characteristic polynomial is $p_k(x)$. The three vertices of $X_1$ have the same neighbors, giving two independent eigenvectors for $0$. The $k-2$ vertices of $Q\setminus\{q_1,q_2,q_3\}$ form a clique and have the same neighbors outside this clique, giving $k-3$ independent eigenvectors for $-1$.

Let $x_2,x_3$ be the unique vertices of $X_2,X_3$, respectively. For either root $\theta$ of $x^2+x-1$, assign coordinates $\theta,-\theta,-1,1$ to $q_2,q_3,x_2,x_3$, respectively, and zero elsewhere. At $q_2$ and $x_2$ the eigenvalue equations are $\theta^2=1-\theta$ and $-\theta=-\theta$. At $q_3$ and $x_3$ they are the corresponding negatives, and at every other vertex the contributions cancel. Thus both roots occur as eigenvalues. Therefore we have accounted for $5+2+(k-3)+2=k+6$ eigenvalues, with multiplicity, which is the order of $H_{k,2}$ for $k\geq3$. This proves~\eqref{eq:hosttwo}.

For $k=2$, the part $Q\setminus\{q_1,q_2,q_3\}$ is empty. Delete its row and column from the previous quotient matrix. The resulting matrix has characteristic polynomial $x^4-x^3-11x^2-2x+12$. The two eigenvectors for $0$ and the two eigenvectors for the roots of $x^2+x-1$ remain, proving~\eqref{eq:hosttwo2}.
\end{proof}
For later root comparisons, note that $p_2(x)=(x+1)(x^4-x^3-11x^2-2x+12)$. Thus the positive roots of $p_2$ are precisely those of the quartic in~\eqref{eq:hosttwo2}, and we may use the same polynomial $p_k$ for every $k\geq2$.

\Needspace{5\baselineskip}
We next consider $H_{k,3}$. Two subdivided branches have the same neighbors outside those branches. Assigning opposite coordinates to the two branches cancels their contributions at the center and at the joined clique. The proof below uses this to construct eigenvectors for $1$ and $-1$.

\begin{lemma}\label{lem:starpoly}
Let $k\geq2$ and $r\geq1$. Then
\begin{equation}\label{eq:starpoly}
 \phi(K_{k-1}\vee S_{r,r};x)=(x+1)^{k+r-3}(x-1)^{r-1}q_{k,r}(x),
\end{equation}
where
\begin{equation}\label{eq:qkr}
 q_{k,r}(x)=x^4+(2-k)x^3+(-2kr-k+r)x^2+(-3kr+k+2r-2)x+(k-1)(r-1)^2.
\end{equation}
\end{lemma}
\begin{proof}
Let $P$ be the vertex set of the $(k-1)$-clique and let $c$ be the center of $S_{r,r}$. Write the subdivided branches as $c u_i v_i$, where $1\leq i\leq r$. The partition with parts $P$, $\{c\}$, $\{u_1,\ldots,u_r\}$, and $\{v_1,\ldots,v_r\}$ is equitable, with quotient matrix
\[
 \begin{pmatrix}
 k-2&1&r&r\\
 k-1&0&r&0\\
 k-1&1&0&1\\
 k-1&0&1&0
 \end{pmatrix}.
\]
Its characteristic polynomial is $q_{k,r}(x)$. The vertices of $P$ have identical neighbors outside $P$, so their coordinate differences give $k-2$ independent eigenvectors for $-1$.

For $1\leq i\leq r-1$, assign coordinates $1,1,-1,-1$ to $u_i,v_i,u_r,v_r$, respectively, and zero elsewhere. Each nonzero coordinate equals the sum of the coordinates at its neighbors. At $c$ and at each vertex of $P$, the contributions from the two branches cancel. Hence this is an eigenvector for $1$. Replacing these four coordinates by $1,-1,-1,1$ gives an eigenvector for $-1$, because the sum of the neighboring coordinates is now the negative of the coordinate at each of the four vertices.

{Therefore we have accounted for $4+(k-2)+(r-1)+(r-1)=k+2r$ eigenvalues, counted with multiplicity,} which is the order of $K_{k-1}\vee S_{r,r}$, proving~\eqref{eq:starpoly}. When $r=1$, there are no pairs of branches to compare, and the two exponents $r-1$ are zero, as required.
\end{proof}

The preceding factorizations give the eigenvalue bounds needed for the sufficient direction of Theorem~\ref{thm:classification}.

\begin{corollary}\label{cor:hosts}
For every $k\geq2$, we have $\lambda_2(H_{k,1})=\lambda_2(H_{k,3})=1$ and $\lambda_2(H_{k,2})<1$. {Also, $\lambda_2(K_k\vee\overline{K_t})=-1$  for $t=0,1$, and $\lambda_2(K_k\vee\overline{K_t})=0$ for $t\geq2$.}
\end{corollary}
\begin{proof}
Equation~\eqref{eq:hostone} gives $\lambda_2(H_{k,1})=1$. For $H_{k,2}$, substitute $x=y+1$ in $p_k$ to obtain
\[
 p_k(y+1)=y^5+(7-k)y^4+(20-11k)y^3+(29-34k)y^2+(15-31k)y-k.
\]
For $k\geq2$, the coefficients of $y^3,y^2,y$ and the constant term are negative. The coefficient of $y^4$ may be positive, zero, or negative. So, $p_k(y+1)$ has exactly one sign change of coefficients, from positive to negative. Therefore it has at most one positive root. Such a root exists because $p_k(1)=-k<0$, whereas $p_k(y+1)$ tends to infinity as $y$ tends to infinity. Thus $p_k$ has exactly one positive root greater than $1$, and $p_k(1)\ne0$ excludes a root at $1$. The roots of the other factors in~\eqref{eq:hosttwo} are less than $1$. The identity for $p_2$ gives the same conclusion from~\eqref{eq:hosttwo2}. Therefore $\lambda_2(H_{k,2})<1$ for every $k\geq2$.

For $H_{k,3}$, use Lemma~\ref{lem:starpoly} with $r=r_k$. Substitution gives
\begin{align*}
 q_{2,9}(x)&=(x-1)(x+4)(x^2-3x-16),\\
 q_{3,8}(y+1)&=y(y^3+3y^2-40y-140),\\
 q_{k,7}(y+1)&=y^4+(6-k)y^3+(19-18k)y^2+(36-53k)y-14
 \quad(k\geq4).
\end{align*}
The quadratic in $q_{2,9}(x)$ has one negative root and one root greater than $1$. The cubic in $q_{3,8}(y+1)$ has one sign change of coefficients, is negative at zero, and has positive leading coefficient, so it has exactly one positive root. The same argument applies to $q_{k,7}(y+1)$, where the coefficients after $y^3$ are negative, and the possible sign of $6-k$ does not introduce another sign change. Hence each $q_{k,r_k}$ has exactly one root greater than $1$. Since $r_k\geq7$, the factor $(x-1)^{r_k-1}$ in~\eqref{eq:starpoly} supplies the eigenvalue $1$, giving $\lambda_2(H_{k,3})=1$.

Finally, let $G=K_k\vee\overline{K_t}$. 
If $t=0$ or $t=1$, then $G$ is a complete graph and $\lambda_2(G)=-1$. If $t\geq2$, then $G=K_k\vee\overline{K_t}=K_{1,1,\ldots,1,t}$ is a complete multipartite graph and $\lambda_2(G)=0$.
\end{proof}

\Needspace{5\baselineskip}
We next calculate the characteristic polynomial that occurs for both graphs in Theorem~\ref{thm:gap}. Its value at $1$ also explains the integer parameters used in their definitions.

\begin{lemma}\label{lem:onearm}
For $k\geq2$ and $t\geq1$,
\begin{equation}\label{eq:onearm}
 \phi(B_k(t,1,0,\ldots,0);x)=x^{t-1}(x+1)^{k-2}f_{k,t}(x),
\end{equation}
where
\begin{equation}\label{eq:fkt}
 f_{k,t}(x)=x^5+(2-k)x^4+\bigl(1-k(t+3)\bigr)x^3-k(t+2)x^2+(2k-1)tx+kt.
\end{equation}
In particular,  $\lambda_2(B_k(t,1,0,\ldots,0))$ is the smaller of the two positive roots of $f_{k,t}$ and is simple. Moreover,  $\lambda_2(B_k(t,1,0,\ldots,0))$ is less than, equal to, or greater than $1$ according as
\begin{equation}\label{eq:fkt-at-one}
 f_{k,t}(1)=(k-1)t-6k+4
\end{equation}
is negative, zero, or positive, respectively.
\end{lemma}
\begin{proof}
The partition with parts $\{q_1\}$, $\{q_2\}$, $Q\setminus\{q_1,q_2\}$, $X_1$, and $X_2$ is equitable, with quotient matrix
\[
 \begin{pmatrix}
 0&1&k-1&0&1\\
 1&0&k-1&t&0\\
 1&1&k-2&t&1\\
 0&1&k-1&0&0\\
 1&0&k-1&0&0
 \end{pmatrix}.
\]
Its characteristic polynomial is $f_{k,t}(x)$. All $t$ vertices of $X_1$ have neighborhood $Q\setminus\{q_1\}$, so their differences give $t-1$ independent eigenvectors for $0$. The $k-1$ vertices of $Q\setminus\{q_1,q_2\}$ form a clique and have identical neighbors outside it, giving $k-2$ independent eigenvectors for $-1$. These vectors have sum zero on the parts of the partition. Together with the five quotient eigenvalues they account for $(t-1)+(k-2)+5=k+t+2$ eigenvalues, proving~\eqref{eq:onearm}.

We next locate the positive roots of $f_{k,t}$. Its leading coefficient is positive, its coefficients of $x^3$ and $x^2$ are negative, and its coefficients of $x$ and the constant term are positive. Since the coefficient $2-k$ of $x^4$ is nonpositive, the nonzero coefficients have exactly two sign changes. Thus there are at most two positive roots, counted with multiplicity. On the other hand,
\[
 f_{k,t}(0)=kt>0,\qquad f_{k,t}(2)=72-48k-(7k+2)t<0,
\]
and $f_{k,t}(x)$ tends to infinity as $x$ tends to infinity. The intermediate value theorem gives a root in $(0,2)$ and a root in $(2,\infty)$. Thus, these are the only positive roots of $f_{k,t}$ and both are simple.

The other factors in~\eqref{eq:onearm} have roots only at $0$ and $-1$. So, the larger positive root of $f_{k,t}$ is the largest eigenvalue of the graph and the smaller one is its second largest eigenvalue. If the smaller root is denoted by $\theta$, then $f_{k,t}$ is positive on $(0,\theta)$ and negative between its two positive roots. The larger root exceeds $2$, so the sign at $1$ is negative, zero, or positive precisely when $\theta<1$, $\theta=1$, or $\theta>1$, respectively. Substituting $x=1$ gives~\eqref{eq:fkt-at-one} and finishes the proof.
\end{proof}

Thus the equation $f_{k,t}(1)=0$ has the real solution $t=6+2/(k-1)$. The largest integer strictly below this value is $\sigma_k$, and the smallest integer strictly above it is $r_k$. The value itself is integral only when $k=2$ or $k=3$.

\section{The structural classification}\label{sec:structure}
In this section, we prove Theorem~\ref{thm:classification}. For each fixed $k\geq2$, we also determine the $k$-trees with $\lambda_2>1$ that are minimal under taking induced subgraphs within the class of $k$-trees for that $k$.

Let $F$ be the $2$-tree with six vertices constructed as follows. Begin with the triangle on $\{u,d,e\}$, add $x$ adjacent to $u,d$, then add $y$ adjacent to $d,e$, and finally add $v$ adjacent to $u,x$. Let $T_0$ be the tree with six vertices obtained by attaching two pendant vertices and a path of length three to a single vertex. Define
\begin{equation}\label{eq:forbidden}
 \begin{aligned}
 F_{1,k}&=K_{k-2}\vee F,& F_{2,k}&=K_{k-1}\vee T_0,
 &F_{3,k}&=K_{k-1}\vee P_6,\\
 F_{4,k}&=B_k(3,2,0,\ldots,0),&
 F_{5,k}&=B_k(4,1,1,0,\ldots,0),
 \end{aligned}
\end{equation}
and, for $k\geq3$, let $F_{6,k}=B_k(3,1,1,1,0,\ldots,0)$. Recall that $E_k=B_k(r_k,1,0,\ldots,0)$, and put
$\Fcal_2=\{F_{1,2},\ldots,F_{5,2},E_2\}$ and
$\Fcal_k=\{F_{1,k},\ldots,F_{6,k},E_k\}$ for $k\geq3$.

We verify that $F_{1,k},F_{2,k}$ and $F_{3,k}$ are $k$-trees.
Let $D$ be a $d$-tree, where $1\leq d\leq k$. Start with the
joining clique $K_{k-d}$ and the initial $d$-clique of $D$. Then they induce $K_k$. Follow a recursive construction
of $D$. Whenever a new vertex is attached to a $d$-clique $C$,
join it also to all vertices of $K_{k-d}$. Its earlier neighbors
then form the $k$-clique $C\cup V(K_{k-d})$. Thus
$K_{k-d}\vee D$ is a $k$-tree. Taking $(d,D)=(2,F)$,
$(1,T_0)$ and $(1,P_6)$ proves that $F_{1,k}$,
$F_{2,k}$ and $F_{3,k}$ are $k$-trees. The remaining members
of $\Fcal_k$ are $k$-trees by the definition of $B_k$.

\Needspace{5\baselineskip}
We first show that these graphs are forbidden. The slightly stronger estimates below will also be used in Section~\ref{sec:gap}.

\begin{lemma}\label{lem:forbidden}
For every $k\geq2$, we have $\lambda_2(F_{i,k})>33/32$ for $1\leq i\leq5$ and $\lambda_2(E_k)>1$. For $k\geq3$, we also have $\lambda_2(F_{6,k})>1.02$.
\end{lemma}
\begin{proof}
For $1\leq i\leq5$, we have $F_{i,2}\ind F_{i,k}$, and for $k\geq3$ we have $F_{6,3}\ind F_{6,k}$. Therefore interlacing reduces the first five bounds to graphs with $k=2$ and the sixth bound to $k=3$. We will handle $E_k$ separately.

We begin with $F_{1,2}=F$, $F_{2,2}=K_1\vee T_0$ and $F_{3,2}=K_1\vee P_6$. Expanding these adjacency determinants, and then substituting $x=y+1$, gives
\begin{align*}
 \phi(F_{1,2};y+1)&=h(y)(y^3+2y^2-5y-9),\\
 \phi(F_{2,2};y+1)&=(y+1)u_2(y),\\
 \phi(F_{3,2};y+1)&=h(y)(y^4+3y^3-5y^2-18y-7),
\end{align*}
where $h(y)=y^3+4y^2+3y-1$ and $u_2(y)=y^6+6y^5+4y^4-34y^3-66y^2-24y+7$.

We have $h(1/32)=-29567/32768<0$ and $h(1)=7>0$. Thus $h$ has a root in $(1/32,1)$, which gives an  eigenvalue in $(33/32,2)$ for both $F_{1,2}$ and $F_{3,2}$. The other factor for $F_{1,2}$ has value $-11$ at $y=1$, and the other factor for $F_{3,2}$ has value $-26$ there. Both have positive leading coefficient, so each also has a root greater than $1$, giving an eigenvalue greater than $2$. So, $F_{1,2}$ and $F_{3,2}$ have at least two eigenvalues greater than $33/32$, i.e., $\lambda_2(F_{1,2})>\frac{33}{32}$ and $\lambda_2(F_{3,2})>\frac{33}{32}$.

For $F_{2,2}$, the values $u_2(1/32)=6640570561/1073741824>0$ and $u_2(1)=-106<0$ give a root in $(1/32,1)$. Since $u_2(1)<0$ and $u_2(y)$ tends to infinity, it has another root greater than $1$. Thus two eigenvalues again lie in $(33/32,2)$ and $(2,\infty)$. Hence, $\lambda_2(F_{2,2})>\frac{33}{32}$.

For $F_{4,2}$, $F_{5,2}$ and $F_{6,3}$ we consider the following equitable quotient matrices. For $F_{4,2}=B_2(3,2,0)$, use the parts $\{q_1\},\{q_2\},\{q_3\},X_1,X_2$. Its quotient matrix is
\[
 Q_4=\begin{pmatrix}
 0&1&1&0&2\\1&0&1&3&0\\1&1&0&3&2\\
 0&1&1&0&0\\1&0&1&0&0
 \end{pmatrix}.
\]
The coordinate differences within $X_1$ and $X_2$ give three eigenvectors for $0$.

For $F_{5,2}=B_2(4,1,1)$, use the parts $\{q_1\},\{q_2,q_3\},X_1,X_2\cup X_3$. For $F_{6,3}=B_3(3,1,1,1)$, use $\{q_1\},\{q_2,q_3,q_4\},X_1,X_2\cup X_3\cup X_4$. The respective quotient matrices are
\[
 Q_5=\begin{pmatrix}0&2&0&2\\1&1&4&1\\0&2&0&0\\1&1&0&0\end{pmatrix},
 \qquad
 Q_6=\begin{pmatrix}0&3&0&3\\1&2&3&2\\0&3&0&0\\1&2&0&0\end{pmatrix}.
\]
In $F_{5,2}$, the four vertices of $X_1$ give three eigenvectors for $0$, and the construction on $q_2,q_3,X_2,X_3$ in Lemma~\ref{lem:hostB} gives the two roots of $x^2+x-1$. In $F_{6,3}$, the three vertices of $X_1$ give two eigenvectors for $0$. For either root $\theta$ of $x^2+x-1$, use coordinates $\theta,-\theta,-1,1$ on $q_i,q_4,x_i,x_4$, where $x_j$ is the vertex of $X_j$ and $i=2,3$. These give two independent eigenvectors for each root. Their coordinates sum to zero on every part, so in each graph these eigenvectors and the quotient eigenvalues account for all eigenvalues.

Expanding $\det((y+1)I-Q_i)$ gives $u_i(y)$ for $i=4,5,6$, where
\begin{align*}
 u_4(y)&=y^5+5y^4-3y^3-41y^2-40y+6,\\
 u_5(y)&=y^4+3y^3-10y^2-27y+1,\\
 u_6(y)&=y^4+2y^3-19y^2-46y+1.
\end{align*}
Therefore
\begin{align*}
 \phi(F_{4,2};y+1)&=(y+1)^3u_4(y),\\
 \phi(F_{5,2};y+1)&=(y+1)^3(y^2+3y+1)u_5(y),\\
 \phi(F_{6,3};y+1)&=(y+1)^2(y^2+3y+1)^2u_6(y).
\end{align*}
Substitution in the polynomials yields
\[
 u_4(1/32)=\frac{158037153}{33554432}>0,\quad
 u_5(1/32)=\frac{153697}{1048576}>0,\quad
 u_6(1/50)=\frac{452601}{6250000}>0,
\]
whereas $u_4(1)=-72$, $u_5(1)=-32$ and $u_6(1)=-61$. Each polynomial thus has a root between its indicated endpoint and $1$, and a  root greater than $1$. After replacing $y$ by $x-1$, $F_{4,2}$ and $F_{5,2}$ have two eigenvalues greater than $33/32$, and $F_{6,3}$ has two eigenvalues greater than $1.02$. Interlacing gives the required bounds for all $k$.

Finally, the definition of $r_k$ gives $r_k>6+2/(k-1)$. Thus $f_{k,r_k}(1)>0$, and Lemma~\ref{lem:onearm} gives $\lambda_2(E_k)>1$.
\end{proof}

The first forbidden graph reduces the possible structure of a $k$-tree to two families.

\begin{lemma}\label{lem:structural}
Let $k\geq2$, and let $G$ be a $k$-tree containing no induced subgraph isomorphic to $F_{1,k}$. Then $G\cong B_k(a_1,\ldots,a_{k+1})$ for some nonnegative integers $a_i$ $(i=1, \cdots, k+1)$, or $G\cong K_{k-1}\vee T$ for some tree $T$.
\end{lemma}
\begin{proof}
Fix $k\geq2$. We use induction on $n=|V(G)|$. If $n=k$, then $G=K_k\cong K_{k-1}\vee K_1$, so the conclusion holds. Suppose that $n>k$, and assume that the conclusion holds for every $k$-tree of order $n-1$ containing no induced $F_{1,k}$. Choose $v$ to be the last vertex in a recursive construction of $G$, and put $G'=G-v$. Then $G'$ is a $k$-tree of order $n-1$, and $N_G(v)$ is a $k$-clique in $G'$. Since $G'$ is an induced subgraph of $G$, it also contains no induced $F_{1,k}$. The induction hypothesis gives one of the following two cases.

\emph{Case 1: $G'\cong B_k(a_1,\ldots,a_{k+1})$.}
Consider the core $Q$ and attachment classes from the definition of $B_k$. If $N_G(v)\subseteq Q$, then $N_G(v)=Q\setminus\{q_j\}$ for some $j$. Adding $v$ increases $a_j$ by one, so $G$ is again a graph of the required $B_k$ form.

Suppose that $N_G(v)$ is not contained in $Q$. Since the vertices outside $Q$ are independent and $N_G(v)$ is a clique, precisely one vertex $x$ of $N_G(v)$ lies outside $Q$. Put $U=N_G(v)\setminus\{x\}$. Then $|U|=k-1$, and there are distinct $d,e\in Q$ such that $Q=U\cup\{d,e\}$, $N_{G'}(x)=U\cup\{d\}$ and $N_G(v)=U\cup\{x\}$.

We claim that every old attachment vertex is adjacent to all of $U$. Otherwise, there are $y\in V(G')\setminus Q$ and $u_0\in U$ such that $y$ is not adjacent to $u_0$. By the definition of an attachment class, $y$ is adjacent to all of $Q\setminus\{u_0\}$. Then the six vertices $u_0,d,e,x,y,v$ induce $F$, where the first three form a triangle, $x$ is adjacent to $u_0,d$, $y$ is adjacent to $d,e$, and $v$ is adjacent to $u_0,x$. Every vertex of $U\setminus\{u_0\}$ is adjacent to all six vertices. Hence these vertices together induce $K_{k-2}\vee F=F_{1,k}$, a contradiction.

It follows that every old attachment vertex is adjacent to all of $U$ and to exactly one of $d,e$. Consider the graph $G-U$. The vertices $d,e$ are joined by an edge. Every old attachment
vertex other than $x$ is a leaf at $d$ or $e$, while $x$ is
adjacent to $d$ and the new leaf $v$. There are no further edges
among these vertices.  Thus $G-U$ is a tree, and every vertex of $U$ is adjacent to every vertex outside $U$. Therefore $G\cong K_{k-1}\vee T$ for the tree $T=G-U$.

\emph{Case 2: $G'\cong K_{k-1}\vee T$.}
Let $P$ be the vertex set of the $(k-1)$-clique. If $P\subseteq N_G(v)$, then $N_G(v)=P\cup\{x\}$ for some $x\in V(T)$. In this case $G\cong K_{k-1}\vee T'$, where $T'$ is obtained from $T$ by adding a leaf at $x$.

Suppose that $P$ is not contained in $N_G(v)$. A clique in a tree has at most two vertices. Since $N_G(v)$ has $k$ vertices, it must be $P_0\cup\{x,y\}$, where $|P_0|=k-2$ and $xy\in E(T)$. Let $p$ be the unique vertex of $P\setminus P_0$.

If every vertex of $T-\{x,y\}$ is adjacent to $x$ or $y$, then $T$ is a \emph{double star}, namely an edge with leaves attached to its endpoints. Take $P\cup\{x,y\}$ as a core. A leaf at $x$ is adjacent to every core vertex except $y$, a leaf at $y$ misses only $x$, and $v$ misses only $p$. All vertices outside this core are pairwise nonadjacent. Hence $G\cong B_k(a_1,\ldots,a_{k+1})$. 

Now we assume that some vertex of $T-\{x,y\}$ is adjacent to neither $x$ nor $y$. After interchanging $x,y$ if necessary, the component of $T-xy$ containing $x$ has a path $xuw$ of length two. The vertices $p,x,y,u,v,w$ induce $F$, where $p,x,y$ form a triangle, $u$ is adjacent to $p,x$, $v$ to $x,y$, and $w$ to $p,u$. Every vertex of $P_0$ is adjacent to all six vertices, again giving an induced $F_{1,k}$. This contradiction completes Case 2 and the induction.
\end{proof}

The next lemma uses $F_{2,k}$ and $F_{3,k}$ to restrict the tree in Lemma~\ref{lem:structural}. Together, the first three forbidden graphs leave only the two families needed in the classification.

\begin{lemma}\label{lem:twofamilies}
Let $k\geq2$, and let $G$ be a $k$-tree containing no induced subgraph isomorphic to $F_{1,k}$, $F_{2,k}$ or $F_{3,k}$. Then $G\cong B_k(a_1,\ldots,a_{k+1})$ or $G\cong K_{k-1}\vee S_{r,s}$ for some nonnegative integers $a_i$ $(i=1, \cdots, k+1)$ and $0\leq s\leq r$.
\end{lemma}
\begin{proof}
By Lemma~\ref{lem:structural}, it is enough to consider $G=K_{k-1}\vee T$ for a tree $T$. If $T$ contains $T_0$ or $P_6$ as an induced subgraph, then together with the joined clique, it would give $F_{2,k}$ or $F_{3,k}$ in $G$. Hence $T$ contains neither $T_0$ nor $P_6$. Since every path in a tree is induced, the absence of $P_6$ implies that the diameter of $T$ is at most four.

If $T$ has one vertex, then $G=K_k=K_{k-1}\vee S_{0,0}$. So, we may assume that $T$ has at least two vertices. If $T$ has diameter at most three, it is a double star, allowing zero leaves at an endpoint. Its central edge together with the joined clique forms a core for a $B_k(a_1, \ldots, a_{k+1})$.

Now suppose that $T$ has diameter four, and choose a diametral path $abcde$. Note that every vertex of $T$ is within distance two from $c$. Otherwise, a vertex farther from $c$ gives a path of length at least five, contrary to the diameter.

Suppose that a neighbor $u$ of $c$ has two distinct neighbors $x,y$ {other than $c$. Note} that at least one of the paths $cba,cde$ avoids $u$, and we denote it by $czw$. The six vertices $x,y,u,c,z,w$ are distinct, and, because $T$ is a tree, they induce precisely $T_0$. This is a contradiction. Thus every neighbor of $c$ has at most one neighbor other than $c$. Each branch at $c$ is consequently a path of length one or two. If $r$ is the number of branches and $s$ the number of branches of length two, then $T\cong S_{r,s}$.
\end{proof}

\Needspace{10\baselineskip}
We now determine the permitted parameters. For two vectors of the same length, write $\aone\preceq\boldsymbol{b}$ if each coordinate of $\aone$ is at most the corresponding coordinate of $\boldsymbol{b}$.

\begin{proposition}\label{prop:parameters}
Let $k\geq2$. For $\boldsymbol{a}=(a_1, \cdots, a_{k+1})$ with integers $a_1\geq\cdots\geq a_{k+1}\geq0$, $\lambda_2(B_k(\aone))\leq1$ if and only if at least one of the following holds:
\begin{equation}\label{eq:parameters}
 a_2=0,\qquad \aone\preceq(2,\ldots,2),\qquad \aone\preceq(3,1,1,0,\ldots,0),\qquad \aone\preceq(r_k-1,1,0,\ldots,0).
\end{equation}
Also, $\lambda_2(K_{k-1}\vee S_{r,s})\leq1$ if and only if $s=0$ or $r\leq r_k$. 
\end{proposition}
\begin{proof}
{Put $G=B_k(\aone)$ and first assume that $\lambda_2(G)\leq1$. By Lemma~\ref{lem:forbidden} and interlacing, $G$ contains no member of $\Fcal_k$ as an induced subgraph. In the selections of attachment vertices below, we also retain all vertices of the core $Q$.} If $a_2=0$, then only $X_1$ can be nonempty. The set $Q\setminus\{q_1\}$ is a $k$-clique, while $\{q_1\}\cup X_1$ is independent and adjacent to every vertex in $Q\setminus \{q_1\}$. Therefore $B_k(\aone)\cong K_k\vee\overline{K_{a_1+1}}$, and we have the first condition.

Suppose next that $a_2\geq2$. If $a_1\geq3$, then selecting three vertices of $X_1$ and two of $X_2$ gives $F_{4,k}$. Hence $a_1\leq2$. Since the parameters are ordered, every parameter is at most $2$, which is the second condition.

It remains to consider $a_2=1$. If $a_3=0$, all later parameters are zero. If $a_1\geq r_k$, then it would give an induced $E_k$ by selecting $r_k$ vertices of $X_1$ and the vertex of $X_2$. Thus $a_1\leq r_k-1$, giving the fourth condition.

If $a_3=1$, then $a_1\geq4$ would similarly give $F_{5,k}$. Hence $a_1\leq3$. When $k=2$, there are no further coordinates, so the third condition follows. When $k\geq3$ and $a_4=0$, it follows in the same way. Finally, if $k\geq3$ and $a_4=1$, having $a_1=3$ would give $F_{6,k}$. Thus $a_1\leq2$, and we return to the second condition. These cases explain all ordered parameter vectors $\aone$.

We now prove sufficiency. {The first condition gives $B_k(\aone)\cong K_k\vee\overline{K_{a_1+1}}$.} {In this case, Corollary~\ref{cor:hosts} gives $\lambda_2(G)<1$.} The second and third conditions give induced subgraphs of $H_{k,1}$ and $H_{k,2}$, respectively. {In both cases, Corollary~\ref{cor:hosts} and interlacing give $\lambda_2(G)\leq1$.} For the fourth condition, if $a_2=0$ we are in the first case. So, we assume that $a_2=1$. Put $a_1=t$ and we will show
\begin{equation}\label{eq:overlap}
 B_k(t,1,0,\ldots,0)\cong K_{k-1}\vee S_{t+1,1}\qquad(t\geq1).
\end{equation}
To verify this, delete $Q\setminus\{q_1,q_2\}$. The remaining graph has the edge $q_1q_2$, the $t$ vertices of $X_1$ as leaves at $q_2$, and the unique vertex of $X_2$ as a leaf at $q_1$. It is $S_{t+1,1}$ with center $q_2$. {Every deleted core vertex is adjacent to every vertex of this tree, proving~\eqref{eq:overlap}. If $t\leq r_k-1$, retain one subdivided branch of $S_{r_k,r_k}$ in full and the internal vertex from each of $t$ other branches, together with the center. These vertices induce $S_{t+1,1}$. Retaining the joined $(k-1)$-clique as well gives $G\ind H_{k,3}$. Hence Corollary~\ref{cor:hosts} and interlacing give $\lambda_2(G)\leq1$.}

Finally, let $G=K_{k-1}\vee S_{r,s}$. If $s=0$, the center and the joined clique form a $k$-clique, and the remaining $r$ vertices are independent. Thus $G\cong K_k\vee\overline{K_r}$ and $\lambda_2(G)\leq1$ by Corollary~\ref{cor:hosts}. Suppose that $s\geq1$. If $r>r_k$, then $S_{r,s}$ contains $S_{r_k+1,1}$ as an induced subgraph. By~\eqref{eq:overlap}, adjoining the $(k-1)$-clique gives an induced $E_k=B_k(r_k, 1, 0, \ldots, 0)$. This is impossible when $G$ avoids $\Fcal_k$, and also when $\lambda_2(G)\leq1$. Hence $r\leq r_k$.

{Conversely, assume that $r\leq r_k$. Choose $r$ branches of $S_{r_k,r_k}$. Keep $s$ of them in full, delete the endpoint from each of the other $r-s$ chosen branches, and discard all unchosen branches. Together with the center, the remaining vertices induce $S_{r,s}$. Thus $G\ind H_{k,3}$, and Corollary~\ref{cor:hosts} and interlacing give $\lambda_2(G)\leq1$.} 
\end{proof}

\noindent\emph{Remark.}
In the necessity argument above, the assumption
$\lambda_2(G)\leq1$ is used only to exclude the members of
$\Fcal_k$. Thus, for either reduced family, the absence of
these graphs as induced subgraphs already forces the stated
parameter conditions.

We now combine the structural reduction with the parameter restrictions to prove Theorem~\ref{thm:classification}.

\begin{proof}[Proof of Theorem~\ref{thm:classification}]
 Let $G$ be a $k$-tree with $\lambda_2(G)\leq1$. By Lemma~\ref{lem:forbidden} and interlacing, $G$ contains no member of $\Fcal_k$. Lemma~\ref{lem:twofamilies} and Proposition~\ref{prop:parameters} imply that $G$ is in (i) or (ii). Conversely, every graph in (i) has second largest eigenvalue less than $1$ by Corollary~\ref{cor:hosts}, and the same corollary and interlacing give $\lambda_2\leq1$ for every graph in (ii).
\end{proof}

The proof also shows that excluding $\Fcal_k$ is sufficient. We next show that none of these forbidden graphs can be omitted. The following property will be used here and again to prove the uniqueness assertion in Theorem~\ref{thm:gap}. For a fixed $k\geq2$, a $k$-tree $G$ is \emph{induced-minimal with $\lambda_2>1$} within the class of $k$-trees if $\lambda_2(G)>1$ and every proper induced subgraph of $G$ that is a $k$-tree for the same $k$ has second largest eigenvalue at most $1$.

\par\begingroup
\begin{remark}\label{rem:extension}
Let $H$ be a proper induced subgraph of a $k$-tree $G$, and suppose that $H$ is also a $k$-tree. Then there is a vertex $v\in V(G)\setminus V(H)$ whose neighborhood in $H$ is a $k$-clique. Thus $G[V(H)\cup\{v\}]$ is again a $k$-tree.

To see this, we first show that a recursive construction of a $k$-tree can start from any prescribed $k$-clique $Q$. We use induction on its order. The assertion holds for $K_k$. Otherwise, let $z$ be the last vertex in a recursive construction and put $S=N_G(z)$. If $z\notin Q$, apply induction to $G-z$ and then add $z$ on $S$. If $z\in Q$, then $Q\setminus\{z\}$ consists of $k-1$ vertices of $S$. Let $w$ be the remaining vertex of $S$. By induction, construct $G-z$ starting from $S$. Instead, start from $Q$, add $w$ adjacent to all of $Q$, and then follow the remaining steps of that construction. The first $k+1$ vertices induce the same complete graph on $S\cup\{z\}$, and no later vertex is adjacent to $z$, since $N_G(z)=S$. This proves the assertion about $Q$.

Now choose a $k$-clique $Q$ in $H$ and construct $G$ starting from $Q$. Let $v$ be the first vertex added outside $H$. Its $k$ earlier neighbors belong to $H$ and form a clique. To show that these are all its neighbors in $H$, put $h=|V(H)|$. Every nonempty induced subgraph of $G$ has a vertex of degree at most $k$, namely a vertex appearing last in this construction. Hence an induced subgraph of order $m\geq k$ has at most $km-\binom{k+1}{2}$ edges. Since $H$ has exactly $kh-\binom{k+1}{2}$ edges, applying this bound to $G[V(H)\cup\{v\}]$ shows that $v$ has at most $k$ neighbors in $H$. Thus its neighborhood in $H$ is precisely the $k$-clique already identified.
\end{remark}
\par\endgroup

We can now identify all $k$-trees that are induced-minimal with $\lambda_2>1$ within the class of $k$-trees for the same $k$. Remark~\ref{rem:extension} allows us to verify this minimality by checking vertex deletions that leave a $k$-tree for the same $k$.

\begin{theorem}\label{thm:minimal}
For a $k$-tree $G$ with $k\geq2$, we have $\lambda_2(G)\leq1$ if and only if $G$ contains no member of $\Fcal_k$ as an induced subgraph. Moreover, the members of $\Fcal_k$ are precisely the induced-minimal $k$-trees with $\lambda_2>1$, and the list is irredundant.
\end{theorem}
\begin{proof}
If $G$ contains a member of $\Fcal_k$, then Lemma~\ref{lem:forbidden} and interlacing give $\lambda_2(G)>1$.  So, $\lambda_2(G)\leq1$ implies that $G$ contains no member of $\Fcal_k$ as an induced subgraph. Conversely, suppose that $G$ contains no member of $\Fcal_k$. Then Lemma~\ref{lem:twofamilies} implies that $G$ is in one of the two reduced families, and the forbidden-subgraph argument in the proof of Proposition~\ref{prop:parameters} gives the permitted parameter ranges. And then, that proposition gives $\lambda_2(G)\leq1$. This proves the equivalence.

Next we prove that each member $F'$ of $\Fcal_k$ is induced-minimal {with $\lambda_2>1$ within} the class of $k$-trees. Let $H$ be a proper induced subgraph of $F'$ that is itself a $k$-tree. {Remark~\ref{rem:extension}} allows us to add a vertex of $F'$ to $H$ while keeping an induced $k$-tree. Repeat this until just one vertex $v$ of $F'$ is missing. Then $H\ind F'-v$, and both $F'$ and $F'-v$ are $k$-trees. Their edge counts differ by $k$, so $d_{F'}(v)=k$. Therefore it is enough to verify $\lambda_2(F'-v)\leq1$ for every vertex of degree $k$ in every listed graph. Then interlacing will give the same bound for $H$.

The possible deletions are listed below. As in the definitions, all omitted coordinates of $B_k$ are zero.
\begin{center}
\begin{tabular}{@{}ll@{}}
\toprule
$F'$&$F'-v$ for $d_{F'}(v)=k$\\
\midrule
$F_{1,k}$&$K_{k-1}\vee P_4$\\
$F_{2,k}$&$K_{k-1}\vee P_5$ or $K_{k-1}\vee S_{3,1}$\\
$F_{3,k}$&$K_{k-1}\vee P_5$\\
$F_{4,k}=B_k(3,2)$&$B_k(2,2)$ or $B_k(3,1)$\\
$F_{5,k}=B_k(4,1,1)$&$B_k(3,1,1)$ or $B_k(4,1)$\\
$F_{6,k}=B_k(3,1,1,1)$, $k\geq3$&$B_k(2,1,1,1)$ or $B_k(3,1,1)$\\
$E_k=B_k(r_k,1)$&$B_k(r_k-1,1)$ or $B_k(r_k)$\\
\bottomrule
\end{tabular}
\end{center}
In $F_{1,k}$, the vertices of degree $k$ are $y,v$ from the definition of $F$. In $F_{2,k}$ and $F_{3,k}$, the vertices of degree $k$ are the leaves of the trees in $T_0$ and $P_6$. For $B_k$ graphs, an attachment vertex has degree $k$. Every core vertex has degree greater than $k$, because at least two attachment classes are nonempty. Thus the table includes every case that needs to be checked.

\par\begingroup
We verify the eigenvalue bound for every graph in the right column by applying Proposition~\ref{prop:parameters}. In the first three rows, $P_4\cong S_{2,1}$ and $P_5\cong S_{2,2}$, so the resulting joins have $(r,s)=(2,1),(2,2)$ or $(3,1)$. Since $r_k\geq7$, all of them satisfy $r\leq r_k$.

For the remaining rows, the parameters of $B_k(2,2)$ and $B_k(2,1,1,1)$ are all at most $2$, while $B_k(3,1,1)=H_{k,2}$. The graphs $B_k(3,1)$, $B_k(4,1)$ and $B_k(r_k-1,1)$ have the form $B_k(t,1,0,\ldots,0)$ with $t\leq r_k-1$, because $r_k-1\geq6$. Finally, $B_k(r_k)$ has only one nonzero attachment parameter, {so it is isomorphic to $K_k\vee\overline{K_{r_k+1}}$.} Thus every graph in the right column satisfies one of the sufficient conditions in Proposition~\ref{prop:parameters}, which gives $\lambda_2(F'-v)\leq1$ in every case. Since $H\ind F'-v$, interlacing also gives $\lambda_2(H)\leq1$. This proves the required minimality.
\par\endgroup

Now let $G$ be any induced-minimal $k$-tree with $\lambda_2(G)>1$. The equivalence already proved implies that $G$ contains a member of $\Fcal_k$. That copy is itself a $k$-tree with second largest eigenvalue greater than $1$, so by the minimality it cannot be a proper subgraph of $G$. Thus $G$ is isomorphic to a member of $\Fcal_k$.

Finally, the graphs in $\Fcal_k$ are pairwise nonisomorphic. Their orders distinguish them except for the pairs $F_{2,k},F_{3,k}$ and $F_{5,k},F_{6,k}$. The first pair has three and two vertices of degree $k$, respectively. In the second pair, $F_{5,k}$ has a vertex of degree $k+2$, whereas $F_{6,k}$ has none. Since every listed graph is induced-minimal, none of them contains another as a proper induced subgraph. Therefore omitting any member would make the characterization false for that graph itself. This proves irredundancy.
\end{proof}

\Needspace{15\baselineskip}
\section{The eigenvalue one and the order bound}\label{sec:boundary}
{In this section, we determine $k$-trees with second largest eigenvalue equal to $1$. We use the two families in Proposition~\ref{prop:parameters} and their eigenvalue equations to find all such graphs and the multiplicity of $1$. We then give the order bound.}

\begin{theorem}\label{thm:boundary}

Let $k\geq2$, and let $G$ be a $k$-tree. Then $\lambda_2(G)=1$ if and only if $G$ belongs to at least one of the following classes. The multiplicity of $1$ is given in each case.
\begin{enumerate}
\item $G\cong B_k(a_1,\ldots,a_{k+1})$, where every $a_i\in\{0,1,2\}$ and at least two parameters equal $2$. In this case, $m_G(1)=|\{i:a_i=2\}|-1$.
\item $G\cong K_{k-1}\vee S_{r,s}$, where $2\leq s\leq r\leq r_k$. In this case, $m_G(1)=s$ if $(k,r)\in\{(2,9),(3,8)\}$, and $m_G(1)=s-1$ otherwise.
\item $G\cong K_{k-1}\vee S_{r_k,1}$, where $k\in\{2,3\}$. In both cases, $m_G(1)=1$.
\end{enumerate}

\end{theorem}
\par\begingroup
\begin{proof}
Suppose first that $G$ belongs to one of (i)--(iii). Proposition~\ref{prop:parameters} gives $\lambda_2(G)\leq1$. Each of these graphs contains $K_{k+1}$, so $\lambda_1(G)\geq k>1$.  So, to prove $\lambda_2(G)=1$, it is enough to show that $1$ is an eigenvalue of $G$.

For (i), choose distinct indices $i,j$ with $a_i=a_j=2$. Assign coordinates $1,-1$ to $q_i,q_j$, respectively, coordinate $-1$ to each vertex of $X_i$, coordinate $1$ to each vertex of $X_j$, and zero elsewhere. At $q_i$ the sum of the neighboring coordinates is $-1+2=1$, and at a vertex of $X_i$ it is $-1$. The sums of the neighboring coordinates at $q_j$ and $X_j$ are $-1$ and $1$, respectively. At every other vertex, the neighboring coordinates sum to zero. Thus the resulting nonzero vector is an eigenvector for $1$.

For (ii), choose two subdivided branches. Assign coordinate $1$ to both noncentral vertices of the first branch, coordinate $-1$ to both noncentral vertices of the second, and zero elsewhere. At each of the four selected vertices, the neighboring coordinates sum to its own coordinate. At the center the sum is $1-1=0$, and at each vertex of the joined clique it is $2-2=0$. At every other vertex, the sum is $0$. So the resulting nonzero vector is an eigenvector for $1$.

For (iii), equation~\eqref{eq:overlap} gives $G\cong B_k(r_k-1,1,0,\ldots,0)$. The pairs $(k,r_k-1)$ are $(2,8)$ and $(3,7)$. In both cases, $f_{k,r_k-1}(1)=0$, so Lemma~\ref{lem:onearm} {shows that $1$ is an eigenvalue of $G$.} Hence every graph in (i)--(iii) has second largest eigenvalue $1$.

Conversely, suppose that $\lambda_2(G)=1$. The graphs $K_k\vee\overline{K_t}$ in Theorem~\ref{thm:classification}(i) have second largest eigenvalue at most $0$ by Corollary~\ref{cor:hosts}, so $G$ is not isomorphic to any of them. Lemma~\ref{lem:forbidden} and interlacing exclude $F_{1,k},F_{2,k},F_{3,k}$ as induced subgraphs. Therefore Lemma~\ref{lem:twofamilies} gives $G\cong B_k(a_1,\ldots,a_{k+1})$ or $G\cong K_{k-1}\vee S_{r,s}$, and Proposition~\ref{prop:parameters} gives the corresponding parameter conditions. We consider these two forms separately. In each case, we identify which of (i)--(iii) occurs and determine the dimension of the entire eigenspace for $1$.

\emph{Case 1: $G=B_k(a_1, \ldots, a_{k+1})$ with $a_1, \ldots, a_{k+1}\geq0$.}
{By relabeling the core vertices and their attachment classes, we may assume that $a_1\geq\cdots\geq a_{k+1}\geq0$.}
Let $\boldsymbol z$ be an eigenvector for $1$, i.e., it satisfies $A(G)\boldsymbol z=\boldsymbol z$, and let $u_i$ be the coordinate of $\boldsymbol z$ at $q_i$ $(i=1,\ldots, k+1)$. Put $S=\displaystyle\sum_{i=1}^{k+1} u_i$ and $y_i=S-u_i$. {The eigenvalue equation at a vertex of $X_i$ gives its coordinate in $\boldsymbol z$ as $y_i$,} since its neighbors are precisely $Q\setminus\{q_i\}$. Define $y_i$ in this way even when $X_i$ is empty, and put $T=\displaystyle\sum_{i=1}^{k+1} a_i y_i$.

{The row of $A(G)\boldsymbol z=\boldsymbol z$ corresponding to $q_i$ gives} $u_i=S-u_i+T-a_i y_i$.  Substituting $u_i=S-y_i$ gives $(a_i-2)y_i=T-S$. Also, by summing the equations $u_i=S-y_i$ over all $i$ we have $\displaystyle\sum_{i=1}^{k+1} y_i=kS$, {i.e.,} $S=\frac{1}{k}\displaystyle\sum_{i=1}^{k+1} y_i$. {Thus, $u_i=S-y_i$ shows that the coordinates $u_i$ are determined by $y_1,\ldots,y_{k+1}$.} 

Summing the equations $(a_i-2)y_i=T-S$ over all $i$ gives $T-2kS=(k+1)(T-S)$, and hence $T=-(k-1)S/k$. Thus $(a_i-2)y_i=T-S$ and $S=\frac{1}{k}\displaystyle\sum_{i=1}^{k+1} y_i$ imply
\begin{equation}\label{eq:one-eigen-equations}
 (a_i-2)y_i+\frac{2k-1}{k^2}\sum_{j=1}^{k+1}y_j=0
 \qquad(i=1,\ldots, k+1).
\end{equation}
We claim that these equations describe the whole eigenspace for $1$. To check this, take any solution $(y_1,\ldots,y_{k+1})$, assign coordinate $\frac1k\displaystyle\sum_{j=1}^{k+1} y_j-y_i(=u_i)$ to $q_i$, and assign coordinate $y_i$ to every vertex of $X_i$. Summing~\eqref{eq:one-eigen-equations} over all $i$ gives $T=-(k-1)S/k$, and substitution verifies every eigenvalue equation. These assignments are inverse to the definition of the $y_i$. Therefore they give a bijective linear map between the solution space of~\eqref{eq:one-eigen-equations} and the eigenspace of $G$ corresponding to $1$.

Suppose first that there exist exactly $h\geq1$ parameters $a_i$ such that they are equal to $2$.  An equation of \eqref{eq:one-eigen-equations} with $a_i=2$ implies {$\displaystyle\sum_{j=1}^{k+1} y_j=0$}.  So,  \eqref{eq:one-eigen-equations} becomes $(a_i-2)y_i=0$, and hence $y_i=0$ if $a_i\ne2$. Thus, {only} $h$ coordinates $y_i$ can be nonzero real numbers subject to the condition that their sum is zero. Hence, one of {the $h$ coordinates} $y_i$ is determined by the others, {i.e.,} the dimension of the eigenspace for $1$ is $h-1$.  Therefore, $m_G(1)=h-1$. 

Since $\lambda_2(G)=1$, this multiplicity is positive and hence $h\geq2$, i.e., there exist at least two parameters $a_i$ such that they are equal to $2$. Then Proposition~\ref{prop:parameters} forces every $a_i\leq2$. Thus $G$ belongs to (i), and the multiplicity is the one in this case.

It remains to consider the case where $a_i\neq2$ for all $i=1,\ldots,k+1$. We show that this gives only the two graphs in (iii). For a nonzero solution of~\eqref{eq:one-eigen-equations}, put $\alpha=\displaystyle\sum_{i=1}^{k+1} y_i$. If $\alpha=0$, then $(a_i-2)y_i=0$ implies $y_i=0$ for all  $i=1,\ldots,k+1$ as $a_i\ne2$. So $\alpha\ne0$. Now, each equation of  \eqref{eq:one-eigen-equations} determines
$y_i=-\frac{2k-1}{k^2(a_i-2)}\alpha$. {Hence} every solution $(y_1, \ldots, y_{k+1})$ of \eqref{eq:one-eigen-equations} is determined by the single scalar $\alpha$. Since a nonzero solution exists, the dimension of the eigenspace for $1$ is one, i.e., $m_G(1)=1$.

Summing the equations $y_i=-\frac{2k-1}{k^2(a_i-2)}\alpha$ over all $i$ and dividing by $\alpha$ gives
\begin{equation}\label{eq:one-eigen-condition}
 1+\frac{2k-1}{k^2}\sum_{i=1}^{k+1}\frac{1}{a_i-2}=0.
\end{equation}
We now use the parameter conditions in Proposition~\ref{prop:parameters}. The case $a_2=0$ is excluded because $B_k(a_1,0,\ldots,0)\cong K_k\vee\overline{K_{a_1+1}}$. If every parameter $a_i$ is $0$ or $1$, let $\ell$ be the number equal to $1$. Then the left side of~\eqref{eq:one-eigen-condition} is $-((k-1)+(2k-1)\ell)/(2k^2)$, which is nonzero. So, this case is not possible. If the parameter vector $(a_1, \ldots, a_{k+1})$ is $(3,1,1,0,\ldots,0)$, then the same left side is $1/(2k)$, again nonzero, i.e., this is also not possible. The only remaining form is $(a_1, \ldots, a_{k+1})=(t,1,0,\ldots,0)$ with $3\leq t\leq r_k-1$. By Lemma~\ref{lem:onearm}, the eigenvalue $1$ occurs only when $(k-1)(t-6)=2$. Its integer solutions are $(k,t)=(2,8)$ and $(3,7)$. Equation~\eqref{eq:overlap} identifies these graphs with those in (iii), and their multiplicity is one, as just proved. This completes Case 1.

{\emph{Case 2: $G=K_{k-1}\vee S_{r,s}$.}
If $s=0$, then $G\cong K_k\vee\overline{K_r}$, which has already been excluded. Thus $s\geq1$, and Proposition~\ref{prop:parameters} gives $r\leq r_k$. Put $p=k-1$.

Let $\boldsymbol z$ be an eigenvector for $1$, i.e., it satisfies $A(G)\boldsymbol z=\boldsymbol z$. We first show that the coordinates of $\boldsymbol z$ on the joined $p$-clique are equal. If $p=1$, then there is only one vertex in the $p$-clique, and hence it is true. So, we may assume that $p\geq2$. Take two vertices of this clique with coordinates $a'$ and $a''$, respectively. Since they have the same neighbors outside those two vertices, subtracting their eigenvalue equations gives $a'-a''=a''-a'$, i.e., $a'=a''$. We denote the common coordinate on the $p$-clique by $a$, and let $c$ be the coordinate of $\boldsymbol z$ at the center of $S_{r,s}$. Let $u_i,v_i$ be the coordinates at the internal vertex and the endpoint of the $i$th subdivided branch, respectively, where $i=1,\ldots,s$.

The endpoint of a subdivided branch is adjacent to the internal vertex and to all $p$ vertices of the joined $p$-clique. The internal vertex is also adjacent to the center. Therefore the corresponding rows of $A(G)\boldsymbol z=\boldsymbol z$ give $v_i=pa+u_i$ and $u_i=pa+c+v_i$. Substituting the first equation into the second gives $c=-2pa$. Each unsubdivided leaf of $S_{r,s}$ is adjacent to the center and to the joined $p$-clique, so its coordinate is $pa+c=-pa$. Thus all coordinates of $\boldsymbol z$ are determined by $a,u_1,\ldots,u_s$.

The row of  $A(G)\boldsymbol z=\boldsymbol z$ corresponding to the center gives $c=pa+\displaystyle\sum_{i=1}^{s}u_i-p(r-s)a$, since there are $r-s$ unsubdivided leaves. Substituting $c=-2pa$ gives
\begin{equation}\label{eq:arm-sum}
 \sum_{i=1}^{s}u_i=(r-s-3)pa.
\end{equation}
A vertex of the joined $p$-clique is adjacent to the other $p-1$ vertices in the $p$-clique and to every vertex of $S_{r,s}$. Therefore, its eigenvalue equation is $a=(p-1)a+c+\displaystyle\sum_{i=1}^{s}(u_i+v_i)-(r-s)pa$. Substituting $c=-2pa$, $v_i=pa+u_i$ and~\eqref{eq:arm-sum} gives $a=((r-7)p-1)a$, i.e.,
\begin{equation}\label{eq:arm-condition}
 \bigl((r-7)p-2\bigr)a=0.
\end{equation}

We claim that~\eqref{eq:arm-sum} and~\eqref{eq:arm-condition} describe the whole eigenspace for $1$. To check this, take any solution $(a,u_1,\ldots,u_s)$. Assign coordinate $a$ to each vertex of the joined $p$-clique, coordinate $-2pa$ to the center, coordinates $u_i,pa+u_i$ to the internal vertex and endpoint of the $i$th subdivided branch, respectively, and coordinate $-pa$ to every unsubdivided leaf. These assignments satisfy the equations at the branch vertices and the unsubdivided leaves. Equation~\eqref{eq:arm-sum} gives the equation at the center, and~\eqref{eq:arm-condition}, together with~\eqref{eq:arm-sum}, gives the equation at each clique vertex. These assignments are inverse to the choice of $a,u_1,\ldots,u_s$ above. Therefore they give a bijective linear map between the solution space of these two equations and the eigenspace of $G$ corresponding to $1$.

Suppose first that $(r-7)p-2\ne0$. Equation~\eqref{eq:arm-condition} gives $a=0$, so~\eqref{eq:arm-sum} becomes $\sum_{i=1}^{s}u_i=0$. Hence one of the $s$ coordinates $u_i$ is determined by the other $s-1$ coordinates, which can be chosen freely. Thus the dimension of the eigenspace for $1$ is $s-1$, i.e., $m_G(1)=s-1$. Since $\lambda_2(G)=1$, this multiplicity is positive and hence $s\geq2$. Therefore $G$ belongs to (ii), with the stated multiplicity.

It remains to consider the case where $(r-7)p-2=0$, i.e., $(r-7)(k-1)=2$. Its integer solutions with $k\geq2$ are $(k,r)=(2,9)$ and $(3,8)$. In this case,~\eqref{eq:arm-condition} is satisfied for every $a$, so only~\eqref{eq:arm-sum} remains. This equation determines $u_s$ once $a,u_1,\ldots,u_{s-1}$ have been chosen. These $s$ variables can be chosen freely, so the dimension of the eigenspace for $1$ is $s$, i.e., $m_G(1)=s$. If $s\geq2$, then $G$ belongs to (ii), with the exceptional multiplicity stated there. If $s=1$, then $r=r_k$ and $k\in\{2,3\}$, so $G$ belongs to (iii) and has multiplicity one. This completes Case 2.

The two cases show that every $k$-tree with second largest eigenvalue $1$ belongs to (i), (ii) or (iii), and give the multiplicity in each case.  This completes the proof.}
\end{proof}

\par\endgroup

{{We next use Theorem~\ref{thm:classification} to determine the largest order of a $k$-tree with $\lambda_2\leq1$ that is not isomorphic to $K_k\vee\overline{K_t}$ for any $t\geq0$.}}

\begin{corollary}\label{cor:order}
Put $M_k=\max\{3k+3,k+2r_k\}$. If $G$ is a $k$-tree with $\lambda_2(G)\leq1$ and $G\not\cong K_k\vee\overline{K_t}$ for every $t\geq0$, then $|V(G)|\leq M_k$. Up to isomorphism, equality holds only for $H_{k,3}$ when $2\leq k\leq5$, and only for $H_{k,1}$ when $k\geq6$. 
\end{corollary}
\par\begingroup
\begin{proof}
By Theorem~\ref{thm:classification}, $G$ is an induced subgraph of $H_{k,1}$, $H_{k,2}$ or $H_{k,3}$. These graphs have orders $3k+3$, $k+6$ and $k+2r_k$, respectively. The second order is smaller than the first. If $2\leq k\leq5$, then $k+2r_k>3k+3$, while $3k+3>k+2r_k$ for $k\geq6$. Thus $|V(G)|\leq M_k$. If equality holds, $G$ has the same order as the unique largest graph among these three. Since it is an induced subgraph of that graph, it must be isomorphic to it.
\end{proof}
\par\endgroup
\begin{remark}
By Corollary \ref{cor:hosts}, the graph $H_{k,3}$ for $2\leq k\leq5$, or the graph $H_{k,1}$ for $k\geq6$, has order $M_k$ and second largest eigenvalue $1$. It is not isomorphic to $K_k\vee\overline{K_t}$ for any $t\geq0$, since these graphs have second largest eigenvalue at most $0$. So, from the above corollary, we know that the least integer $N_k$ such that every $k$-tree of order $n\geq N_k$ with $\lambda_2\leq1$ is isomorphic to $K_k\vee\overline{K_{n-k}}$ is $N_k=M_k+1$. In particular, $N_2=21$, $N_3=20$, $N_4=19$, $N_5=20$, and $N_k=3k+4$ for $k\geq6$.
\end{remark}

\Needspace{12\baselineskip}
\section{The second largest eigenvalues nearest to one}\label{sec:gap}
{In this section, we prove Theorem~\ref{thm:gap}. We first give bounds for $\lambda_2(L_k)$ and $\lambda_2(E_k)$. We then compare the characteristic polynomials of the graphs obtained in Section~\ref{sec:structure}, and determine the graphs attaining the two values.}

{The following bounds will be used in the comparisons and in the proof of the equality cases.}

\begin{lemma}\label{lem:rootbounds}
For every $k\geq2$, $0.97<\lambda_2(L_k)<1<\lambda_2(E_k)<1.03$, and both $\lambda_2(L_k)$ and $\lambda_2(E_k)$ are simple. In particular,
$1.01<\lambda_2(E_2),\lambda_2(E_3)<1.02$, $\lambda_2(E_4)<1.01$ and $\lambda_2(E_k)<\xi$ for $k\geq4$, where $\xi$ is the unique zero in $(1,1.03)$ of the polynomial
\begin{equation}\label{eq:D}
 D(x)=x^4+10x^3+9x^2-14x-7.
\end{equation}
\end{lemma}
\par\begingroup
\begin{proof}
By Lemma~\ref{lem:onearm}, $\lambda_2(L_k)$ and $\lambda_2(E_k)$ are simple and are the smaller positive roots of $f_{k,\sigma_k}$ and $f_{k,r_k}$, respectively. The other positive root of each polynomial is greater than $2$. Also, $\sigma_k<6+2/(k-1)<r_k$, so the same lemma gives $\lambda_2(L_k)<1<\lambda_2(E_k)$.

{We first prove the lower bound for $\lambda_2(L_k)$. At $x=0.97$, the values of $f_{2,7}$ and $f_{3,6}$ are about $0.9517$ and $0.7623$, respectively. So, both are positive.}  By~\eqref{eq:fkt},
\begin{equation}\label{eq:f6}
 f_{k,6}(x)=x\bigl((x^2+x)^2-6\bigr)-k(x-1)T(x),
\end{equation}
where $T(x)=x^3+10x^2+18x+6$.  {At $x=0.97$, the coefficient of $k$ in~\eqref{eq:f6} is positive. Together with $f_{3,6}(0.97)>0$, this gives $f_{k,6}(0.97)>0$ for every $k\geq3$.} Each of these polynomials is negative between its two positive roots. Since its larger positive root exceeds $2$, its smaller positive root must be greater than $0.97$. Thus $\lambda_2(L_k)>0.97$.

Now, we consider $\lambda_2(E_k)$. For $k=2,3,4$, substitution gives
\[
 \begin{array}{c|ccc}
 &f_{2,9}&f_{3,8}&f_{4,7}\\\hline
 x=1.01&+&+&-\\
 x=1.02&-&-&-
 \end{array}
\]
where the entries are the exact signs of the indicated values. All three polynomials are positive at $1$. Their smaller positive roots therefore lie in $(1.01,1.02)$ for $k=2,3$, and in $(1,1.01)$ for $k=4$. This proves the additional bounds.

It remains to give the upper bound for all $k\geq4$. Again by~\eqref{eq:fkt},
\begin{equation}\label{eq:f7}
 f_{k,7}(x)=x\bigl((x^2+x)^2-7\bigr)-kD(x).
\end{equation}
For $x\geq1$, we have $D'(x)=4x^3+30x^2+18x-14\geq38>0$, i.e., $D(x)$ is strictly increasing on $x\geq1$. Also, $D(1)=-1$ and $D(1.03)>0$. {Thus $D$ has a unique root $\xi$ in $(1,1.03)$.} Since $(x^2+x)^2<7$ for $1\leq x\leq1.03$, we obtain $f_{k,7}(1)=k-3>0$ and $f_{k,7}(\xi)<0$. Hence $f_{k,7}$ has a root in $(1,\xi)$. It is the smaller positive root, since the other is greater than $2$. Therefore $1<\lambda_2(E_k)<\xi<1.03$.
\end{proof}
\par\endgroup
\begin{remark}
Exact substitution gives $D(1.025537352055)<0<D(1.025537352056)$. Since $D$ is strictly increasing there, these signs give  $1.025537352055<\xi<1.025537352056$. So, $\xi\approx1.025537352055$.
\end{remark}

\Needspace{6\baselineskip}
{We first prove the bound below $1$ and determine its equality case.}

\begin{proposition}\label{prop:lower}
For every $k\geq2$, every $k$-tree $G$ with $\lambda_2(G)<1$ satisfies $\lambda_2(G)\leq\lambda_2(L_k)$, with equality if and only if $G\cong L_k$.
\end{proposition}
\par\begingroup
\begin{proof}
Let $G$ be a $k$-tree with $\lambda_2(G)<1$. By Lemma~\ref{lem:forbidden} and interlacing, $G$ contains no member of $\Fcal_k$. {Therefore Lemma~\ref{lem:twofamilies} gives $G\cong B_k(a_1,\ldots,a_{k+1})$ or $G\cong K_{k-1}\vee S_{r,s}$. Since $\lambda_2(G)<1$, Proposition~\ref{prop:parameters} gives the corresponding parameter conditions in either case.} We first determine which graphs remain after excluding the equality cases in Theorem~\ref{thm:boundary}.

Suppose that $G=B_k(a_1, \ldots, a_{k+1})$, with $a_1\geq\cdots\geq a_{k+1}$. {If $a_2=0$, then $G\cong K_k\vee\overline{K_{a_1+1}}$, as shown in Proposition~\ref{prop:parameters}.} If all parameters are at most $2$, then by Theorem~\ref{thm:boundary}(i), at most one can equal $2$. In this case $G\ind B_k(2,1,\ldots,1)$. The third parameter condition gives $G\ind H_{k,2}$. In the remaining case, $G=B_k(a,1,0,\ldots,0)$ with $1\leq a\leq r_k-1$. Theorem~\ref{thm:boundary}(iii) excludes $a=8$ when $k=2$ and $a=7$ when $k=3$, so $a\leq\sigma_k$.

Now suppose that $G=K_{k-1}\vee S_{r,s}$. Theorem~\ref{thm:boundary}(ii) gives $s\leq1$. {If $s=0$, then $G\cong K_k\vee\overline{K_r}$. If $s=1$ and $r=1$, then $G\cong K_{k-1}\vee P_3\cong K_k\vee\overline{K_2}$.} For $s=1$ and $r\geq2$, equation~\eqref{eq:overlap} gives $G\cong B_k(r-1,1,0,\ldots,0)$, and the preceding argument applies.

{The graphs $K_k\vee\overline{K_t}$ with $t\geq0$ have second largest eigenvalue at most $0$, by Corollary~\ref{cor:hosts}, whereas $\lambda_2(L_k)>0$ by Lemma~\ref{lem:rootbounds}. Thus it remains to compare $H_{k,2}$, $B_k(2,1,\ldots,1)$ and the graphs $B_k(a,1,0,\ldots,0)$ with $1\leq a\leq\sigma_k$.}

For $B_k(2,1,\ldots,1)$, the partition with parts $\{q_1\}$, $Q\setminus\{q_1\}$, $X_1$ and $X_2\cup\cdots\cup X_{k+1}$ is equitable. Its quotient matrix is
\[
 \begin{pmatrix}
 0&k&0&k\\1&k-1&2&k-1\\0&k&0&0\\1&k-1&0&0
 \end{pmatrix},
\]
with characteristic polynomial
\begin{equation}\label{eq:dk}
 d_k(x)=x^4+(1-k)x^3-(k+1)^2x^2-k(k-1)x+2k^2.
\end{equation}
Assigning $1,-1$ to the two vertices of $X_1$ and zero elsewhere gives an eigenvector for $0$. For $i\geq2$, let $x_i$ be the vertex of $X_i$. For each root $\theta$ of $x^2+x-1$ and each $i=2,\ldots,k$, assign coordinates $\theta,-\theta,-1,1$ to $q_i,q_{k+1},x_i,x_{k+1}$ and zero elsewhere. These are eigenvectors for $\theta$, as in the proof of Lemma~\ref{lem:hostB}. {For a fixed root $\theta$, the vector indexed by $i$ has coordinate $\theta\ne0$ at $q_i$, whereas every other vector in this collection has coordinate $0$ there. If a linear combination of these vectors is zero, its coordinate at $q_i$ forces the coefficient of the $i$th vector to be zero. This holds for every $i=2,\ldots,k$, so the $k-1$ vectors are linearly independent.} Their coordinate sums on every part are zero. Together with the four quotient eigenvalues, these give all $2k+3$ eigenvalues. Hence
\begin{equation}\label{eq:smallones}
 \phi(B_k(2,1,\ldots,1);x)=x(x^2+x-1)^{k-1}d_k(x).
\end{equation}

The nonzero coefficients of both $p_k$ and $d_k$ have two sign changes. Also, both polynomials are positive at $0$, whereas $p_k(1)=-k<0$ and $d_k(1)=1-2k<0$. Since their leading coefficients are positive, each has a root in $(0,1)$ and another in $(1,\infty)$. Descartes' rule shows that these are the only positive roots of $p_k$ and $d_k$, and both are simple. {For $H_{k,2}$, formulas~\eqref{eq:hosttwo}--\eqref{eq:hosttwo2} show that the positive eigenvalues are the two positive roots of $p_k$ and $(\sqrt5-1)/2$. Note that for $k=2$, the extra factor $x+1$ in $p_2$ does not affect the positive roots. Thus $\lambda_2(H_{k,2})$ is the larger of the root of $p_k$ in $(0,1)$ and $(\sqrt5-1)/2$. Similarly, formula~\eqref{eq:smallones} shows that $\lambda_2(B_k(2,1,\ldots,1))$ is the larger of the root of $d_k$ in $(0,1)$ and $(\sqrt5-1)/2$.}

For $1\leq a\leq\sigma_k$, put $\theta_a=\lambda_2(B_k(a,1,0,\ldots,0))$. By Lemma~\ref{lem:onearm}, $\theta_a$ is the unique root of $f_{k,a}$ in $(0,1)$. Put $u_k=1-1/(4k)$. Then $(\sqrt5-1)/2<u_k<1$. Substitution gives
\begin{align*}
 d_k(u_k)&=-\frac{k^4(320k-624)+k^2(400k-132)+20k-1}{256k^4}<0,\\
 f_{k,6}(u_k)&=\frac{k^4(6912k-5184)+k^2(1808k-308)+28k-1}{1024k^5}>0.
\end{align*}
Thus the root of $d_k$ in $(0,1)$ is less than $u_k$, while $u_k<\theta_6$. It follows that $\lambda_2(B_k(2,1,\ldots,1))<\theta_6$.

For $H_{k,2}$, we have
\[
 f_{k,6}(x)-p_k(x)
 =x\bigl\{(k-2)+(1-x)((2k+1)x+3k+2)\bigr\}>0
 \qquad(0<x<1).
\]
Let $\zeta$ be the root of $p_k$ in $(0,1)$. Then $f_{k,6}(\zeta)>0$. Since $\theta_6$ is the unique root of $f_{k,6}$ in $(0,1)$ and the other positive root of $f_{k,6}$ is greater than $2$, this implies that $\zeta<\theta_6$. Also, $(\sqrt5-1)/2<u_k<\theta_6$. Therefore $\lambda_2(H_{k,2})<\theta_6$.

Finally, let $1\leq a<b\leq\sigma_k$. For $0<x<1$,
\[
 f_{k,b}(x)-f_{k,a}(x)
 =(b-a)\bigl\{k(1+x(1-x)(x+2))-x\bigr\}>0.
\]
Note that the expression in braces is the sum of $k-x>0$ and $kx(1-x)(x+2)>0$. Taking $x=\theta_b$ gives $f_{k,a}(\theta_b)<0$, so $\theta_a<\theta_b$. Thus the largest value among $\theta_1, \ldots, \theta_{\sigma_k}$ is $\theta_{\sigma_k}$. {Since $\sigma_k\geq6$, we have $\theta_6\leq\theta_{\sigma_k}=\lambda_2(L_k)$. We have already proved that $\lambda_2(H_{k,2})<\theta_6$ and $\lambda_2(B_k(2,1,\ldots,1))<\theta_6$. If $G\ind H_{k,2}$, then interlacing gives $\lambda_2(G)\leq\lambda_2(H_{k,2})<\theta_{\sigma_k}$. If $G\ind B_k(2,1,\ldots,1)$, then it similarly gives $\lambda_2(G)\leq\lambda_2(B_k(2,1,\ldots,1))<\theta_{\sigma_k}$. The graphs $K_k\vee\overline{K_t}$ also have second largest eigenvalue strictly less than $\theta_{\sigma_k}$. Thus equality can occur only in the remaining family $B_k(a,1,0,\ldots,0)$, where strict increase in $a$ forces $a=\sigma_k$.} Hence $\lambda_2(G)\leq\lambda_2(L_k)$, with equality if and only if $G\cong L_k$.
\end{proof}
\par\endgroup

{We now consider a $k$-tree whose second largest eigenvalue is greater than $1$. By Theorem~\ref{thm:minimal}, it contains a member of $\Fcal_k$ as an induced subgraph. Thus it is enough to compare the second largest eigenvalues of these forbidden graphs.}

\begin{proposition}\label{prop:upper}
For every $k\geq2$, every $k$-tree $G$ with $\lambda_2(G)>1$ satisfies $\lambda_2(G)\geq\lambda_2(E_k)$. If equality holds, then $E_k\ind G$.
\end{proposition}
\par\begingroup
\begin{proof}
Let $G$ be a $k$-tree with $\lambda_2(G)>1$. By Theorem~\ref{thm:minimal}, $G$ contains a member of $\Fcal_k$ as an induced subgraph. For $1\leq i\leq5$, Lemmas~\ref{lem:forbidden} and~\ref{lem:rootbounds} give $\lambda_2(F_{i,k})>33/32>1.03>\lambda_2(E_k)$. For $k=2$, this shows all the required comparisons. For $k\geq3$, it remains to compare $F_{6,k}$ with $E_k$.

We first assume $k\geq4$ and consider the equitable partition of $F_{6,k}$ with parts $\{q_1\}$, $\{q_2,q_3,q_4\}$, $Q\setminus\{q_1,q_2,q_3,q_4\}$, $X_1$ and $X_2\cup X_3\cup X_4$. Its quotient matrix is
\[
 \begin{pmatrix}
 0&3&k-3&0&3\\1&2&k-3&3&2\\1&3&k-4&3&3\\
 0&3&k-3&0&0\\1&2&k-3&0&0
 \end{pmatrix},
\]
with characteristic polynomial
\begin{equation}\label{eq:hk}
 h_k(x)=x^5+(2-k)x^4+(3-8k)x^3+(2-9k)x^2+(10k-9)x+9k.
\end{equation}
The three vertices of $X_1$ have identical neighborhoods, giving two independent eigenvectors for $0$. The $k-3$ vertices of $Q\setminus\{q_1,q_2,q_3,q_4\}$ have identical neighbors outside their clique, giving $k-4$ independent eigenvectors for $-1$. For each root $\theta$ of $x^2+x-1$, assign coordinates $\theta,-\theta,-1,1$ to $q_i,q_4,x_i,x_4$, where $x_j$ is the vertex of $X_j$ and $i=2,3$, and zero elsewhere. As in Lemma~\ref{lem:hostB}, these are two independent eigenvectors for $\theta$. Their coordinates sum to zero on every part. Counting these vectors and the quotient eigenvalues gives all $k+7$ eigenvalues, and therefore
\begin{equation}\label{eq:forbidden-six-poly}
 \phi(F_{6,k};x)=x^2(x+1)^{k-4}(x^2+x-1)^2h_k(x).
\end{equation}
For $k=3$, the third part $Q\setminus\{q_1, q_2, q_3, q_4\}$ is empty. The quotient matrix $Q_6$ in Lemma~\ref{lem:forbidden} gives
\[
 \phi(F_{6,3};x)=x^2(x^2+x-1)^2\frac{h_3(x)}{x+1},
\]
where  $h_3(x)=(x+1)(x^4-2x^3-19x^2-6x+27)$. Thus, for every $k\geq3$, the roots of $h_k$ greater than $1$ are precisely the eigenvalues of $F_{6,k}$ greater than $1$. 

Substituting $x=y+1$ gives
\[
 h_k(y+1)=y^5+(7-k)y^4+(21-12k)y^3+(33-39k)y^2+(17-36k)y+k-1.
\]
For $k\geq3$, the coefficients of $y^3,y^2,y$ are negative, and the leading and constant coefficients are positive. Whatever the sign of $7-k$, the nonzero coefficients have exactly two sign changes. Hence $h_k$ has at most two roots greater than $1$, counted with multiplicity. On the other hand, $h_k(1)=k-1>0$ and $h_k(2)=78-87k<0$, while $h_k(x)$ tends to infinity. So, there is one root in $(1,2)$ and another in $(2,\infty)$. Hence, these are the only two roots greater than $1$, and both are simple. The larger one is the largest eigenvalue of $F_{6,k}$ and is at least $k$, since the core induces $K_{k+1}$.

Put $\gamma=\lambda_2(E_k)$. For $k\geq4$,
\begin{equation}\label{eq:uppercompare}
 h_k(x)-f_{k,7}(x)=2(x^2+x-1)\bigl((k+1)x-k\bigr).
\end{equation}
For $k=3$, the corresponding difference is $h_3(x)-f_{3,8}(x)=(x-1)(11x^2+16x-3)$. Both expressions are positive for $x>1$. Since $f_{k,r_k}(\gamma)=0$, we have $h_k(\gamma)>0$. The polynomial $h_k$ is negative between its two roots greater than $1$, and its larger root is at least $k>\gamma$. Therefore $\gamma$ is less than its smaller root, i.e., $\lambda_2(E_k)<\lambda_2(F_{6,k})$.

We have shown that $E_k$ has smaller second largest eigenvalue than every other member of $\Fcal_k$. So, interlacing gives $\lambda_2(G)\geq\lambda_2(E_k)$. If equality holds, no other member of $\Fcal_k$ can occur in $G$, and hence $E_k\ind G$.
\end{proof}
\par\endgroup

{Proposition~\ref{prop:lower} gives the equality case below $1$. To determine the equality case above $1$, we first consider the eigenvectors for $\lambda_2(E_k)$.}

\begin{lemma}\label{lem:vectors}
Let $k\geq2$. Then every eigenvector for $\lambda_2(E_k)$ has all coordinates nonzero, and the sum of its coordinates on every $k$-clique of $E_k$ is nonzero.
\end{lemma}
\par\begingroup
\begin{proof}
Put $\theta=\lambda_2(E_k)$, $t=r_k$ and $w=\theta(\theta+1)$. By Lemma~\ref{lem:rootbounds}, $1<\theta<1.03$, so $w>2$ and $w^2<6<t$. Let $\boldsymbol z$ be a nonzero eigenvector for $\theta$. The coordinates on $X_1$ are equal, since its vertices have the same neighbors and $\theta\ne0$. The coordinates on $Q\setminus\{q_1,q_2\}$ are also equal.

Let $a,b,c,e,d$ be the coordinates at $q_1,q_2$, on $Q\setminus\{q_1,q_2\}$, on $X_1$, and at the vertex of $X_2$, respectively. The eigenvalue equations give
\[
 \begin{gathered}
 \theta a=b+(k-1)c+d,\qquad \theta b=a+(k-1)c+te,\\
 \theta c=a+b+(k-2)c+te+d,\qquad
 \theta e=b+(k-1)c,\qquad \theta d=a+(k-1)c.
 \end{gathered}
\]
We first show that $e\ne0$. Suppose that $e=0$. The equation at $X_1$ then gives $b=-(k-1)c$. Subtracting the equation for $c$ from the equation for $a$ gives $a=c$. Subtracting the equation for $c$ from the equation for $b$ gives $d=k(\theta+1)c$, while the equation at $X_2$ gives $\theta d=kc$. Hence $k(w-1)c=0$. Since $w>2$, we obtain $c=0$, and then every coordinate is zero. This contradicts that $\boldsymbol z$ is a nonzero eigenvector for $\theta$.

Multiplying $\boldsymbol z$ by a nonzero scalar, we may assume that $e=1$. Subtracting the equation for $c$ from the first two equations gives $(\theta+1)(a-c)=-t$ and $(\theta+1)(b-c)=-d$. Together with the equations at $X_1$ and $X_2$, these imply $k(\theta+1)c=w+d$ and $(w-1)d=w-t$. Therefore
\begin{equation}\label{eq:eigenvector}
 d=\frac{w-t}{w-1},\qquad
 c=\frac{w^2-t}{k(\theta+1)(w-1)},\qquad
 a=c-\frac{t}{\theta+1},\qquad b=\theta-(k-1)c.
\end{equation}
Since $w^2<t$ and $w>2$, we have $a,c,d<0<b$. Since the coordinates on $X_1$ are $e=1$,  every coordinate of $\boldsymbol z$ is nonzero.

It remains to check the sum of its coordinates on each $k$-clique. Suppose first that the clique is contained in $Q$. For $Q\setminus\{q_1\}$ and $Q\setminus\{q_2\}$, the sums are $\theta$ and $\theta d$, respectively, by the equations at the attachment vertices. For $Q\setminus\{q_j\}$ with $j\geq3$, the sum is $a+b+(k-2)c=\theta+a-c=(w-t)/(\theta+1)$. None of these sums is zero.

Now suppose that the clique contains an attachment vertex. Then there is exactly one such vertex, since the attachment vertices are independent. If it belongs to $X_2$, then the other vertices are in $Q\setminus\{q_2\}$, and all its coordinates are negative. If it belongs to $X_1$ and the clique contains $q_2$, then the sum is $1+b+(k-2)c=1+\theta-c>0$.

The remaining case consists of one vertex of $X_1$ and all $k-1$ vertices of $Q\setminus\{q_1,q_2\}$. Subtracting the equation at $X_2$ from the equation at $q_2$ gives $b=d+t/\theta$. Thus the sum is
\[
 1+(k-1)c=1+\theta-b
 =\frac{(\theta^2-1)(w-t)}{\theta(w-1)}<0.\qedhere
\]
\end{proof}
\par\endgroup

{We now complete the proof of Theorem~\ref{thm:gap}. For the equality case above $1$, we use Remark~\ref{rem:extension}, Lemma~\ref{lem:vectors} and interlacing.}

\par\begingroup
\begin{proof}[Proof of Theorem~\ref{thm:gap}]
{Proposition~\ref{prop:lower} shows that $\lambda_2(G)\leq\lambda_2(L_k)$ for every $k$-tree $G$ with $\lambda_2(G)<1$. Proposition~\ref{prop:upper} shows that $\lambda_2(G)\geq\lambda_2(E_k)$ for every $k$-tree $G$ with $\lambda_2(G)>1$.} Lemma~\ref{lem:rootbounds} gives $\lambda_2(L_k)<1<\lambda_2(E_k)$. Since both $L_k$ and $E_k$ are $k$-trees, they attain the two bounds. Hence $\beta_k$ and $\gamma_k$ exist, and $\beta_k=\lambda_2(L_k)$ and $\gamma_k=\lambda_2(E_k)$. Proposition~\ref{prop:lower} also gives the equality case for $\beta_k$.

It remains to prove the equality case for $\gamma_k$. Suppose that $\lambda_2(G)=\gamma_k$. By Proposition~\ref{prop:upper}, $G$ contains an induced copy of $E_k$. If this copy $E_k$ is proper, then Remark~\ref{rem:extension} gives a vertex $v\in V(G)$  outside the copy of $E_k$ whose neighborhood in $E_k$ is a $k$-clique. Put $H=G[V(E_k)\cup\{v\}]$.

{Let $\boldsymbol z$ be a nonzero eigenvector of $E_k$ corresponding to $\gamma_k$. By Lemma~\ref{lem:vectors}, $\sum_{u\in N_H(v)}z_u\ne0$. Interlacing gives $\lambda_2(H)\geq\lambda_2(E_k)=\gamma_k$. Suppose that equality holds. The equality condition of interlacing stated in Section~\ref{sec:spectral} gives a nonzero eigenvector $\boldsymbol y$ of $E_k$ corresponding to $\gamma_k$ whose extension by zero at $v$ is an eigenvector of $H$. Since $\gamma_k$ is simple by Lemma~\ref{lem:rootbounds}, we have $\boldsymbol y=c\boldsymbol z$ for some $c\ne0$. The eigenvalue equation at $v$ is then $0=\sum_{u\in N_H(v)}y_u=c\sum_{u\in N_H(v)}z_u$, a contradiction. Therefore $\lambda_2(H)>\gamma_k$.} Since $H\ind G$, interlacing gives $\lambda_2(G)>\gamma_k$, a contradiction. Therefore $G\cong E_k$. Conversely, $E_k$ attains $\gamma_k$ by definition. This completes the proof.
\end{proof}
\par\endgroup

{Finally, we determine how $\beta_k$ and $\gamma_k$ vary with $k$.}

\begin{corollary}\label{cor:limits}
\begin{itemize}
\item[(i)] The sequence $(\beta_k)_{k\geq3}$ is strictly increasing and tends to $1$. Moreover, $k(1-\beta_k)$ tends to $2/35$. 
\item[(ii) ]The sequence $(\gamma_k)_{k\geq4}$ is strictly increasing and tends to the unique zero $\xi$ in $(1,1.03)$ of the polynomial $D$ defined in~\eqref{eq:D}. Moreover, $\gamma_4<\min\{\gamma_2,\gamma_3\}$.

\item[(iii)]For  $\mathcal S=\{\lambda_2(G):G\text{ is a $k$-tree for some }k\geq2\}$, 
\[
 \sup\bigl(\mathcal S\cap(-\infty,1)\bigr)=1,
 \qquad
 \min\bigl(\mathcal S\cap(1,\infty)\bigr)=\gamma_4>1.
\]
\end{itemize}
\end{corollary}
\par\begingroup
\begin{proof}
We first consider $\beta_k$ for $k\geq3$. Since $\sigma_k=6$ for $k\geq3$, we have $f_{k,6}(\beta_k)=0$. By~\eqref{eq:f6}, $f_{k+1,6}(\beta_k)=(1-\beta_k)T(\beta_k)>0$. The smaller positive root of $f_{k+1,6}$ is $\beta_{k+1}<1$, and the larger is greater than $2$. The polynomial is negative between these roots. {Since $\beta_k<1$, $f_{k+1,6}(\beta_k)>0$ implies} that $\beta_k<\beta_{k+1}$, and hence the sequence is increasing and bounded above by $1$. Thus, the sequence has a limit $\ell\leq1$. Equation~\eqref{eq:f6} gives
\[
 k(1-\beta_k)
 =\frac{\beta_k\bigl(6-(\beta_k^2+\beta_k)^2\bigr)}{T(\beta_k)}.
\]
Since $0.97<\beta_k<1$ and $T$ is positive on $[0.97,1]$, the right side is bounded. If $\ell<1$, the left side would tend to infinity, a contradiction. Thus $\ell=1$. Then the right side tends to $2/T(1)=2/35$.

Now we consider $\gamma_k$ for $k\geq4$. Lemma~\ref{lem:rootbounds} gives $1<\gamma_k<\xi$. Since $D$ is increasing on $x\geq1$ and $D(\xi)=0$, we have $D(\gamma_k)<0$. Equation~\eqref{eq:f7} gives $f_{k+1,7}(\gamma_k)=-D(\gamma_k)>0$. As above, the smaller positive root of $f_{k+1,7}$ is $\gamma_{k+1}$ and the larger exceeds $2$, whereas $\gamma_k<1.03$. Hence $\gamma_k<\gamma_{k+1}$. The sequence is bounded above by $\xi$ and therefore has a limit. {Also, from equation~\eqref{eq:f7}, we have}
\[
 D(\gamma_k)=\frac{\gamma_k\bigl((\gamma_k^2+\gamma_k)^2-7\bigr)}{k}.
\]
The numerator is bounded, so the right side tends to zero. The limit of $(\gamma_k)_{k\geq4}$ lies in $[1,\xi]$, and $\xi$ is the only zero of $D$ in that interval. Therefore $\gamma_k$ tends to $\xi$.

Since $\beta_k$ tends to $1$, we know $ \sup\bigl(\mathcal S\cap(-\infty,1)\bigr)=1$. Now we determine the minimum above $1$. Since the sequence $(\gamma_k)_{k\geq4}$ is strictly increasing, $\gamma_4<\gamma_k$ for every $k\geq5$. {Using $f_{4,7}(\gamma_4)=0$ and subtracting the relevant polynomials gives}
\[
 \begin{aligned}
 f_{2,9}(\gamma_4)&=2(\gamma_4-1)(\gamma_4^3+9\gamma_4^2+16\gamma_4+5)>0,\\
 f_{3,8}(\gamma_4)&=1+(\gamma_4-1)(\gamma_4^3+8\gamma_4^2+14\gamma_4+5)>0.
 \end{aligned}
\]
Thus $f_{k,r_k}(\gamma_4)>0$ for $k=2,3$. Since $\gamma_4<2$, the same root comparison as above gives $\gamma_4<\gamma_k$ for $k=2,3$. Thus, $\min\bigl(\mathcal S\cap(1,\infty)\bigr)=\gamma_4>1$.
\end{proof}
\par\endgroup

\section{Uniform forbidden induced subgraphs}\label{sec:uniform}
For each fixed $k\geq2$, Theorem~\ref{thm:minimal} determines the $k$-trees with $\lambda_2>1$ that are induced-minimal within the class of $k$-trees for the same $k$. In this section, we consider $k$-trees $G$ with $\lambda_2(G)>1$ such that every proper induced subgraph $H$ of order at least two satisfies $\lambda_2(H)\leq1$, where $H$ need not be connected or be a $k$-tree. The original graph $G$ is still required to be a $k$-tree for some $k\geq1$.

{We also replace the lists $\Fcal_k$, whose graphs depend on $k$, by one list of nine graphs that applies to every $k\geq1$. Theorem~\ref{thm:uniform} shows that this list is irredundant and that its members are precisely the $k$-trees satisfying the stronger minimality condition. Each listed graph has at most $13$ vertices. Therefore every $k$-tree $G$ with $\lambda_2(G)>1$ contains a connected induced subgraph $H$ of order at most $13$ with $\lambda_2(H)>1$, independently of $k$ and $|V(G)|$. Corollary~\ref{cor:thirteen} shows that this order bound is best possible.}

{We need the tree case in order to include $k=1$. For nonnegative integers $a,b$, let $D_{a,b}$ be the double star obtained from an edge by attaching $a$ leaves to one endpoint and $b$ leaves to the other. Part (i) of the following proposition is the forbidden-subtree characterization stated in the remark following Neumaier's Theorem~4.6~\cite{Neumaier}. Part (ii) follows from his Theorems~2.6, 4.4 and 4.6. Part (iii) is an elementary consequence of (i), characteristic polynomials and interlacing. We include the calculations for (iii) to compare the tree case with Theorem~\ref{thm:gap} and to calculate the eigenvalues used in the proof of Theorem~\ref{thm:uniform}. These tree results are preliminary to the common forbidden-subgraph list.}

\begin{proposition}\label{prop:trees}
Let $T$ be a tree of order at least two. Then the following hold.
\begin{enumerate}
\item $\lambda_2(T)\leq1$ if and only if $T$ contains none of $P_6,T_0,D_{3,2}$ as an induced subgraph.
\item $\lambda_2(T)=1$ if and only if $T\cong D_{2,2}$ or $T\cong S_{r,s}$ with $2\leq s\leq r$.
\item  $\displaystyle\sup_{T\in\mathcal{T_-}}\lambda_2(T)=1$ and $\displaystyle\min_{T\in\mathcal{T_+}}\lambda_2(T)=\sqrt{3-\sqrt3}$, where $\mathcal{T_-}=\{T\mid T\text{ is a tree}, |V(T)|\geq 2, \lambda_2(T)<1\}$ and $\mathcal{T_+}=\{T\mid T\text{ is a tree}, |V(T)|\geq 2, \lambda_2(T)>1\}$. The supremum is not attained, and the minimum is attained
only by $D_{3,2}$ up to isomorphism.

\end{enumerate}
\end{proposition}
\par\begingroup
\begin{proof}
We use Neumaier's results cited above for (i) and (ii), and prove (iii).

First let $a,b\geq1$. Partition the vertices of $D_{a,b}$ into the two singleton centers and the two sets of leaves. This partition is equitable, with quotient matrix
\[
 \begin{pmatrix}0&1&a&0\\1&0&0&b\\1&0&0&0\\0&1&0&0\end{pmatrix}.
\]
Its characteristic polynomial is $x^4-(a+b+1)x^2+ab$. At each center, fix one of its leaves and pair it with each of the other leaves at that center. For each pair, assign coordinates $1,-1$ to the two leaves and zero elsewhere. Since the two leaves have the same neighbor, this gives an eigenvector for $0$. The two centers give $a-1$ and $b-1$ independent vectors, respectively. If a center has only one leaf, it contributes no such vector. Together with the four quotient eigenvalues, these vectors account for all $a+b+2$ eigenvalues, counted with multiplicity. Therefore

\begin{equation}\label{eq:double-star}
 \phi(D_{a,b};x)=x^{a+b-2}\bigl(x^4-(a+b+1)x^2+ab\bigr).
\end{equation}
The two roots of $z^2-(a+b+1)z+ab$ are positive and distinct, since their sum and product are positive and the discriminant is $(a-b)^2+2a+2b+1>0$. Thus the quartic has two distinct positive roots, and $\lambda_2(D_{a,b})$ is the smaller one.

We also calculate the values for the three trees in (i). Removing a short leaf of $T_0$ gives $P_5$, while removing this leaf and its neighbor gives $K_1\cup P_3$. Hence~\eqref{eq:leaf-recurrence} gives $\phi(T_0;x)=x\phi(P_5;x)-x\phi(P_3;x)=x^2(x^4-5x^2+5)$. Together with the path spectrum~\cite[Section~1.4.4]{BrouwerHaemers} and~\eqref{eq:double-star}, this gives
\begin{equation}\label{eq:tree-values}
 \lambda_2(P_6)=2\cos\frac{2\pi}{7},\qquad
 \lambda_2(T_0)=\sqrt{\frac{5-\sqrt5}{2}},\qquad
 \lambda_2(D_{3,2})=\sqrt{3-\sqrt3}.
\end{equation}
All three values are greater than $1$. Moreover, $2\pi/7<3\pi/10$ gives $2\cos(2\pi/7)>2\cos(3\pi/10)=\sqrt{(5-\sqrt5)/2}$. Also, $\sqrt3>(1+\sqrt5)/2$ gives $(5-\sqrt5)/2>3-\sqrt3$. Thus the smallest of the three values is $\lambda_2(D_{3,2})$.

For the values below $1$, consider $S_{r,1}\cong D_{r-1,1}$ with $r\geq2$. Equation~\eqref{eq:double-star} gives
\begin{equation}\label{eq:tree-lower}
 \lambda_2(S_{r,1})^2=\frac{r+1-\sqrt{(r-1)^2+4}}{2}<1.
\end{equation}
These second largest eigenvalues are positive. Rationalizing the difference gives
\[
 1-\lambda_2(S_{r,1})^2
 =\frac{2}{\sqrt{(r-1)^2+4}+r-1}.
\]
The right-hand side tends to zero as $r$ goes to infinity. Since $0<\lambda_2(S_{r,1})<1$, it follows that $\lambda_2(S_{r,1})$ tends to $1$ from below. Therefore, for any tree $T$ with $\lambda_2(T)<1$, we can  choose $r$ sufficiently
large that $\lambda_2(T)<\lambda_2(S_{r,1})<1$. Hence the supremum of the values below $1$ is $1$, and none of these values is a maximum.

Now suppose that $\lambda_2(T)>1$. By (i), $T$ contains $P_6$, $T_0$ or $D_{3,2}$ as an induced subgraph. The comparison in~\eqref{eq:tree-values} and interlacing give $\lambda_2(T)\geq\lambda_2(D_{3,2})$. So, the graph $D_{3,2}$ attains this bound.

It remains to determine the equality case. Suppose that $\lambda_2(T)=\lambda_2(D_{3,2})$. The strict comparisons above exclude $P_6$ and $T_0$, so $T$ contains an induced copy of $D_{3,2}$. If this copy is proper, connectedness gives a vertex $w$ outside it adjacent to it. The vertex $w$ has exactly one neighbor in the copy, since two neighbors would give a cycle.

By~\eqref{eq:double-star}, $\lambda_2(D_{3,2})$ is simple and positive. Let $\boldsymbol z$ be a corresponding nonzero eigenvector. We will show that every coordinate of $\boldsymbol z$ is nonzero. If a center has coordinate zero, the eigenvalue equations at its leaves force all their coordinates to be zero. The equation at that center then forces the coordinate at the other center to be zero, and hence all the other leaf coordinates are zero. If a leaf has coordinate zero, its
eigenvalue equation forces the coordinate at its center to be zero, and the same argument applies. In either case the entire vector would be zero, a contradiction.

Apply the strict interlacing argument in Section~\ref{sec:spectral} to this copy of $D_{3,2}$ and the new vertex $w$. The eigenvalue $\lambda_2(D_{3,2})$ is simple, and the coordinate of $\boldsymbol z$ at the unique neighbor of $w$ is nonzero. Therefore $\lambda_2\bigl(T[V(D_{3,2})\cup\{w\}]\bigr)>\lambda_2(D_{3,2})$. Since this graph is an induced subgraph of $T$, interlacing gives $\lambda_2(T)>\lambda_2(D_{3,2})$, a contradiction. Thus $T\cong D_{3,2}$, completing the proof of (iii).
\end{proof}
\par\endgroup

{The bound $\lambda_2(G)\geq\gamma_4$ also holds for $k=1$ whenever $\lambda_2(G)>1$. Indeed, Proposition~\ref{prop:trees} gives $\lambda_2(G)\geq\sqrt{3-\sqrt3}>1.03>\gamma_4$ in this case. For $k\geq2$, use Theorem~\ref{thm:gap} and Corollary~\ref{cor:limits}. Therefore $E_4$ is the unique graph attaining the smallest value above $1$ among all $k$-trees with $k\geq1$.}

We now return to arbitrary $k\geq1$. The minimality condition in Theorem~\ref{thm:minimal} does not imply that every vertex deletion has second largest eigenvalue at most $1$. For example, $F_{2,2}=K_1\vee T_0$ is induced-minimal with
$\lambda_2>1$ within the class of $2$-trees by that theorem. Let $v$ be its joined vertex. Then $F_{2,2}-v\cong T_0$, which is not a $2$-tree and has second largest eigenvalue greater than $1$ by~\eqref{eq:tree-values}. Thus the minimality condition in Theorem~\ref{thm:minimal} holds for $F_{2,2}$, but the condition on all vertex deletions does not.

{For a graph $G$ of order at least three, the stronger condition is}
\begin{equation}\label{eq:vertex-deletion}
 \lambda_2(G)>1,\qquad
 \lambda_2(G-v)\leq1\quad\text{for every }v\in V(G).
\end{equation}
{If every vertex deletion satisfies $\lambda_2\leq1$, let $H$ be any proper induced subgraph of order at least two and choose a vertex $v$ outside it. Then $H\ind G-v$, so $\lambda_2(H)\leq1$ by interlacing. Conversely, the condition on all proper induced subgraphs includes every vertex deletion. Thus~\eqref{eq:vertex-deletion} is equivalent to requiring $\lambda_2(G)>1$ and $\lambda_2(H)\leq1$ for every proper induced subgraph $H$ of order at least two. Unlike Theorem~\ref{thm:minimal}, we do not require $H$ to be a $k$-tree.} Put
\begin{equation}\label{eq:uniform-list}
 \mathcal O=\{P_6,T_0,D_{3,2},F_{1,2},F_{5,2},F_{6,3},E_2,E_3,E_4\}.
\end{equation}
{Each member of $\mathcal O$ is a $k$-tree for some $k\in\{1,2,3,4\}$. The same list will be used for every $k$, without requiring a forbidden induced subgraph to be a $k$-tree for that value of $k$.}

\begin{theorem}\label{thm:uniform}
Let $G$ be a $k$-tree of order at least two, where $k\geq1$. Then $\lambda_2(G)\leq1$ if and only if $G$ contains no member of $\mathcal O$ as an induced subgraph. {The list is irredundant within the class of all $k$-trees with $k\geq1$.} Moreover, a $k$-tree of order at least three satisfies~\eqref{eq:vertex-deletion} if and only if it is isomorphic to a member of $\mathcal O$. In particular, no $k$-tree with $k\geq5$ satisfies \eqref{eq:vertex-deletion}.
\end{theorem}
\par\begingroup
\begin{proof}
If $G$ contains a member of $\mathcal O$, then Lemma~\ref{lem:forbidden}, equation~\eqref{eq:tree-values}, and interlacing give $\lambda_2(G)>1$. Hence $\lambda_2(G)\leq1$ implies that no member of $\mathcal O$ occurs as an induced subgraph.

Conversely, suppose that $\lambda_2(G)>1$. If $k=1$, then Proposition~\ref{prop:trees} gives an induced $P_6,T_0$ or $D_{3,2}$. Now assume that $k\geq2$. By Theorem~\ref{thm:minimal}, $G$ contains a member of $\Fcal_k$. We will show that each member of $\Fcal_k$ contains a graph in $\mathcal O$ as an induced subgraph.

{First consider $F_{1,k},F_{2,k},F_{3,k}$. Deleting their joined cliques leaves $F=F_{1,2}$, $T_0$ and $P_6$, respectively, all of which belong to $\mathcal O$.}

{For $F_{4,k}=B_k(3,2,0,\ldots,0)$, consider $q_1,q_2$ as two central vertices and all five attachment vertices as the leaves. Then the resulting graph is $D_{3,2}$. For $F_{5,k}$, $q_1,q_2,q_3$ and all attachment vertices give $F_{5,2}$. Similarly, for $F_{6,k}$ with $k\geq3$, we obtain $F_{6,3}$ from $q_1,q_2,q_3,q_4$ and all attachment vertices.}

{Finally, consider $E_k$. The graphs $E_2$ and $E_3$ are already in $\mathcal O$. If $k\geq4$, then $r_k=7$. From $q_1,q_2$, any three other core vertices, and all eight attachment vertices we obtain $B_4(7,1,0,0,0)=E_4$. Thus every member of $\Fcal_k$ contains a member of $\mathcal O$ as an induced subgraph, and hence so does $G$. This proves the equivalence.}

{The equivalence does not yet prove minimality under arbitrary vertex deletion. We now verify~\eqref{eq:vertex-deletion} for every member of $\mathcal O$, including deletions that do not leave a $k$-tree for the same $k$. We first consider the three trees.} Deleting a leaf from $P_6$ gives $P_5\cong S_{2,2}$. Deleting a leaf from $T_0$ gives either $P_5$ or $S_{3,1}\cong D_{2,1}$. Deleting a leaf from $D_{3,2}$ gives $D_{2,2}$ or $D_{3,1}$. Proposition~\ref{prop:trees} and~\eqref{eq:double-star} give second largest eigenvalue at most $1$ for each of these graphs.

A nonleaf deletion from $P_6$ gives $K_1\cup P_4$ or $P_2\cup P_3$. A nonleaf deletion from $T_0$ gives $2K_1\cup P_3$, $P_2\cup P_3$ or $K_1\cup K_{1,3}$. Deleting a center of $D_{3,2}$ gives isolated vertices together with a star. Note that the spectrum of each union is the union of the spectra of its components. Isolated vertices and $P_2$ have no eigenvalue greater than $1$, and each of $P_3,P_4,K_{1,3}$ and a star has at most one eigenvalue greater than $1$. Thus every nonleaf deletion has at most one eigenvalue greater than 1, giving second largest eigenvalue at most $1$.

Now consider $F_{1,2},F_{5,2},F_{6,3},E_2,E_3,E_4$. The five graphs of the form $B_k(a_1,\ldots,a_{k+1})$ are $F_{5,2}=B_2(4,1,1)$, $F_{6,3}=B_3(3,1,1,1)$, $E_2=B_2(9,1,0)$, $E_3=B_3(8,1,0,0)$, and $E_4=B_4(7,1,0,0,0)$. For each of these five graphs, Theorem~\ref{thm:minimal} verifies every attachment-vertex deletion for the corresponding value of $k$. Theorem~\ref{thm:minimal} also verifies deletion of $y$ or $v$ from $F_{1,2}=F$. So, it remains to check the core vertices in the five $B_k$ graphs and the vertices $u,d,e,x$ in $F$.

For a core deletion from $H=B_k(a_1,\ldots,a_{k+1})$, let the deleted vertex be $q_j$. Put $m=\sum_i a_i$ and $I_+=\{i:a_i>0\}$. Partition $H-q_j$ into the $k$ singleton core vertices and the nonempty sets $X_\ell$, and denote the quotient matrix by $Q_j$. The remaining core is a $k$-clique. A core vertex $q_i$ has $a_\ell$ neighbors in $X_\ell$ if $i\ne\ell$, and none otherwise. A vertex of $X_\ell$ is adjacent to every remaining core vertex other than $q_\ell$, when $q_\ell$ is present. There are no edges between attachment vertices. These neighbor counts determine $Q_j$.

Within a set $X_\ell$, fix one vertex and assign coordinates $1,-1$ to it and one other vertex, with zero coordinates elsewhere. This gives $a_\ell-1$ independent eigenvectors for $0$. Summing over all nonempty attachment classes gives $m-|I_+|$ such eigenvectors, each with coordinate sum zero on every part. Together with the $k+|I_+|$ eigenvalues of $Q_j$, they account for all $k+m$ eigenvalues. Hence
\[
 \phi(H-q_j;x)=x^{m-|I_+|}\det(xI-Q_j).
\]

For example, let $H=F_{5,2}=B_2(4,1,1)$ and delete $q_1$. Use the parts $\{q_2\}$, $\{q_3\}$, $X_1$, $X_2$ and $X_3$, in this order, as in the construction above. The quotient matrix is \[
Q_1=\begin{pmatrix}
0&1&4&0&1\\
1&0&4&1&0\\
1&1&0&0&0\\
0&1&0&0&0\\
1&0&0&0&0
\end{pmatrix},
\]
and expanding its characteristic polynomial gives $\det(xI-Q_1)=x(x^2+x-1)(x^2-x-9)$. Note that in this case, $m=6$ and $|I_+|=3$, so the factor contributed by coordinate differences within the attachment classes is $x^{m-|I_+|}=x^3$. Hence
$\phi(F_{5,2}-q_1;x)=x^4(x^2+x-1)(x^2-x-9)$.

The deletions from $F$ can be computed directly from its definition. For example, $F-u$ has edge set $\{de,dx,dy,ey,xv\}$, while $F-e$ has edge set $\{ud,ux,dx,dy,uv,xv\}$.

For each core vertex $q_j$ of $F_{5,2},F_{6,3},E_2,E_3$ or $E_4$, calculate $x^{m-|I_+|}\det(xI-Q_j)$ using the quotient matrix $Q_j$ defined above. For the deletions of $u,d,e,x$ from $F$, calculate the adjacency determinants directly.
After substituting $x=y+1$, these factorizations give the table below, including the case $F_{5,2}-q_1$ calculated above. In each row,

\begin{equation}\label{eq:deletion-factor}
 \phi(H-v;y+1)=y^{\alpha}(y+1)^{\beta}(y+2)^{\delta}
 (y^2+3y+1)^{\varepsilon}g_{H,v}(y).
\end{equation}
The vertex labels of $F$ are those in its definition. The permutation $(u\ d)(e\ x)(y\ v)$ is an automorphism of $F$, so deletion of $u$ or $d$ gives isomorphic graphs, as does deletion of $e$ or $x$. In a $B_k$ graph, core vertices with equal attachment parameters can be exchanged together with their attachment classes. Thus it is enough to check one vertex in each row.

\Needspace{9\baselineskip}
\begingroup
\small
\setlength{\tabcolsep}{3pt}
\renewcommand{\arraystretch}{1.15}
\begin{longtable}{@{}llcl@{}}
\toprule
$H$ & $v$ & $(\alpha,\beta,\delta,\varepsilon)$ & $g_{H,v}(y)$\\
\midrule
\endfirsthead
\toprule
$H$ & $v$ & $(\alpha,\beta,\delta,\varepsilon)$ & $g_{H,v}(y)$\\
\midrule
\endhead
\bottomrule
\endfoot
$F_{1,2}$ & $u,d$ & $(1,0,1,0)$ & $y^{3} + 3 y^{2} - y - 5$ \\
 & $e,x$ & $(0,0,1,0)$ & $y^{4} + 3 y^{3} - 2 y^{2} - 8 y - 2$ \\
\addlinespace[3pt]
$F_{5,2}$ & $q_1$ & $(0,4,0,1)$ & $y^{2} + y - 9$ \\
 & $q_2,q_3$ & $(1,4,0,0)$ & $y^{3} + 4 y^{2} - 2 y - 14$ \\
\addlinespace[3pt]
$F_{6,3}$ & $q_1$ & $(0,3,0,2)$ & $y^{2} - 14$ \\
 & $q_2,q_3,q_4$ & $(1,3,0,1)$ & $y^{3} + 3 y^{2} - 11 y - 32$ \\
\addlinespace[3pt]
$E_2$ & $q_1$ & $(0,8,0,0)$ & $y^{4} + 4 y^{3} - 14 y^{2} - 54 y - 28$ \\
 & $q_2$ & $(0,8,1,0)$ & $y^{3} + 2 y^{2} - 10 y - 2$ \\
 & $q_3$ & $(0,8,0,0)$ & $y^{4} + 4 y^{3} - 5 y^{2} - 18 y - 1$ \\
\addlinespace[3pt]
$E_3$ & $q_1$ & $(0,7,1,0)$ & $y^{4} + 3 y^{3} - 25 y^{2} - 79 y - 36$ \\
 & $q_2$ & $(0,7,2,0)$ & $y^{3} + y^{2} - 20 y - 4$ \\
 & $q_3,q_4$ & $(1,7,0,0)$ & $y^{4} + 5 y^{3} - 11 y^{2} - 73 y - 74$ \\
\addlinespace[3pt]
$E_4$ & $q_1$ & $(0,6,2,0)$ & $y^{4} + 2 y^{3} - 34 y^{2} - 98 y - 42$ \\
 & $q_2$ & $(0,6,3,0)$ & $y^{3} - 28 y - 6$ \\
 & $q_3,q_4,q_5$ & $(1,6,1,0)$ & $y^{4} + 4 y^{3} - 23 y^{2} - 110 y - 105$ \\
\end{longtable}
\endgroup

Every factor in~\eqref{eq:deletion-factor} other than $g_{H,v}(y)$ has no positive root in $y$. Since $x=y+1$, the corresponding adjacency eigenvalues are at most $1$.
Each $g_{H,v}$ has exactly one sign change in its nonzero coefficients, i.e., there exists at most one positive root, counted with multiplicity. Such a root exists because $g_{H,v}(0)<0$ and $g_{H,v}(y)$ tends to infinity. Therefore $\phi(H-v;y+1)$ has exactly one positive root. Thus, $H-v$ has exactly one eigenvalue greater than $1$, i.e., $\lambda_2(H-v)\leq1$.

The table covers every remaining vertex of the six graphs. Together with the degree-$k$ deletions and the tree deletions already checked, it proves~\eqref{eq:vertex-deletion} for every member of $\mathcal O$.

Now let a $k$-tree $G$ satisfy~\eqref{eq:vertex-deletion}. By the equivalence proved above, $G$ contains an induced copy of some $H\in\mathcal O$. If this copy is proper, choose $v\in V(G)\setminus V(H)$. Then $H\ind G-v$, so interlacing gives $\lambda_2(G-v)\geq\lambda_2(H)>1$, a contradiction. Thus $G\cong H$.

Finally, the three trees in $\mathcal O$ are $1$-trees; $F_{1,2},F_{5,2},E_2$ are $2$-trees; $F_{6,3},E_3$ are $3$-trees; and $E_4$ is a $4$-tree. Their orders and clique numbers distinguish them except for $P_6$ and $T_0$, whose degree sequences differ. Hence they are pairwise nonisomorphic. No member can contain another as a proper induced subgraph, since deleting a vertex outside that copy would contradict~\eqref{eq:vertex-deletion}. Thus omitting any member makes the characterization false for that graph itself, proving irredundancy.

Every $k$-tree other than the initial graph $K_k$ has clique number $k+1$. All graphs in $\mathcal O$ have clique number at most five, so none is a $k$-tree with $k\geq5$. For $k\geq5$, the initial graph $K_k$ has second largest eigenvalue $-1$ and does not satisfy~\eqref{eq:vertex-deletion}. Therefore no $k$-tree with $k\geq5$ satisfies~\eqref{eq:vertex-deletion}.
\end{proof}
\par\endgroup

\Needspace{11\baselineskip}
{The orders of the graphs in $\Fcal_k$ depend on $k$, whereas every member of $\mathcal O$ has order at most $13$. Thus the previous theorem gives the following order bound for every $k\geq1$.}

\begin{corollary}\label{cor:thirteen}
Let $G$ be a $k$-tree with $k\geq1$ and $\lambda_2(G)>1$. Then $G$ has a connected induced subgraph $H$ with $2\leq|V(H)|\leq13$ and $\lambda_2(H)>1$, where the constant $13$ is best possible.
\end{corollary}
\par\begingroup
\begin{proof}
By Theorem~\ref{thm:uniform}, $G$ contains a member $H$ of $\mathcal O$ as an induced subgraph. Every member of $\mathcal O$ is connected, has at most $13$ vertices, and has second largest eigenvalue greater than $1$. Thus $H$ has all the required properties.

To show that $13$ is best possible, consider $E_4$, which has $13$ vertices. By Theorem~\ref{thm:uniform}, every vertex deletion from $E_4$ has second largest eigenvalue at most $1$. If $H$ is a proper induced subgraph of $E_4$ of order at least two, choose a vertex $v$ outside it. Since $H\ind E_4-v$, interlacing gives $\lambda_2(H)\leq\lambda_2(E_4-v)\leq1$. Therefore no induced subgraph of $E_4$ with at most $12$ vertices satisfies the required inequality. Hence the bound cannot be reduced.
\end{proof}
\par\endgroup
\begin{remark}
The above corollary also says the following: Let $G$ be a $k$-tree of order at least two, where $k\geq1$. Then $\lambda_2(G)\leq1$ if and only if every connected induced subgraph $H$ of $G$ with $2\leq|V(H)|\leq13$ satisfies $\lambda_2(H)\leq1$. 
\end{remark}

\section*{Acknowledgements}
This research is supported in part by the NSFC (No.~12171154) and the Program of China Scholarship Council (No.~202506740024). J. Park was supported by the National Research Foundation of Korea (NRF) grant funded by the Korea government (MSIT) (RS-2024-00356153).

\end{document}